\documentclass[preprint,3p]{elsarticle}

\usepackage{afterpage}
\usepackage{amssymb}
\usepackage{amsmath}
\usepackage{bm}
\usepackage{cancel}
\usepackage{centernot}
\usepackage{empheq}
\usepackage{float}
\usepackage{lineno}
\usepackage{physics}
\usepackage{siunitx}
\usepackage{soul}
\usepackage{subfig}
\usepackage{subfloat}
\usepackage{tikz}
\usepackage{todonotes}
\usepackage[normalem]{ulem}
\usetikzlibrary{arrows,positioning,shapes.geometric}
\DeclareSIUnit\Molar{M}

\begin{document}

\begin{frontmatter}



\title{A mathematical model for irreversible damage of the collagen scaffold in the myocardium}


\author[a]{Irena Radišić\corref{cor1}}
\ead{irena.radisic@polimi.it}
\author[a]{Francesco Regazzoni}
\ead{francesco.regazzoni@polimi.it}
\author[a]{Luca Dede'}
\ead{luca.dede@polimi.it}
\author[b,c]{Alfio Quarteroni}
\ead{alfio.quarteroni@polimi.it}

\cortext[cor1]{Corresponding author}
\affiliation[a]{organization={MOX -- Department of Mathematics, Politecnico di Milano},
            addressline={Piazza Leonardo da Vinci, 32}, 
            city={Milano},
            postcode={20133},
            country={Italy}}
\affiliation[b]{organization={MOX -- Department of Mathematics (Professor Emeritus), Politecnico di Milano},
            addressline={Piazza Leonardo da Vinci, 32}, 
            city={Milano},
            postcode={20133},
            country={Italy}}
\affiliation[c]{organization={Mathematics Institute (Professor Emeritus), École Polytechnique Fédérale de Lausanne},
            postcode={CH-1015},
            country={Switzerland}}

\begin{abstract}
We propose a dissipative, irreversible damage model for anisotropic media in large deformations to address the damage process of the myocardium following a cardiac infarction. We model damage to the collagen in the cleavage planes  resulting from an increased load to the passive tissue. Starting from variational principles, we derive a quasi-static differential model governing the irreversible evolution of the damage. We apply our model to two test cases. First, a rectangular passive slab geometry to which an indenter-like load is applied. Then, we couple our model with a model for the numerical simulation of left ventricular electromechanics. The computational results obtained with our novel mathematical model show that a redistribution of the load within a left ventricle occurring following a cardiac infarction can produce the damage to the collagen scaffold in the infarcted area consistent with experimental results. Our model and results support the evidence that passive load-dependent dissipative phenomena are implicated in post-infarction remodelling.
\end{abstract}


\begin{keyword}
Damage modelling \sep Nonlinear mechanics \sep Cardiac mechanics \sep Dissipative processes \sep Cardiac electromechanics


\end{keyword}

\end{frontmatter}



\section{Introduction}
\label{sec:introduction}
The heart is an organ responsible for the supply of oxygenated blood throughout the body. It is composed of four chambers, where fluid passes from one chamber into another through the cardiac valves. The largest chamber, and the one with the thickest wall, is the left ventricle, where re-oxygenated blood exits through the aorta in order to be supplied to the rest of the body. By tissue composition, the heart is a muscle primarily composed of myocytes, cardiac muscle cells, that are interconnected by a stiff collagen scaffold. The space between myocytes is filled with fluid and other interstitial components~\cite{katz_physiology_2011}. The myocardium contracts due to the microscopic cellular electrical activity, which propagates across the myocardium and activates the cells individually. The change in cellular polarisation initiates a series of events which causes a microscopic sliding of filaments resulting in a macroscopic muscle contraction. Other than this active mechanical behaviour, the myocardium exhibits an anisotropic passive material behaviour, that is consequence of its structure. Cardiomyocytes are connected in parallel by endomysial collagen strands in order to form mesostructural laminar sheetlets of about four cardiomyocytes each across the thickness. These sheetlets are further separated by cleavage planes made of perimysial collagen fibres~\cite{pope_three-dimensional_2008}. Albeit myocytes are not perfectly parallelly aligned, it is possible to identify a local principal mean direction, which rotates along the myocardium~\cite{wilson_myocardial_2022}. Interstitial fluid and other extracellular components offer an isotropic matrix for the cardiomyocytes and collagen scaffold to be immersed into.
\par
Oxygenated blood is supplied to the myocardium by the coronary arteries. An occlusion of the coronary circulation causes a myocardial infarction~\cite{katz_physiology_2011,hall_guyton_2011}, yielding local loss of contractility and structural changes both in the following days and in the following weeks. Prolonged oxygen starvation of the affected tissue can cause cellular death. The functional and morphological changes to the left ventricle are referred to as postinfarction left ventricular remodelling~\cite{sutton_left_2000}.
\par
The specific biochemical signal driving the remodelling process depends on the stage and single process of the postinfarction remodelling. The affected tissue, which has no longer a viable oxygen supply, loses its ability to actively contract. Due to the local lack of active contraction, the tissue poses less resistance to the fluid pressure in the systolic phase. Consequently, the infarcted region suffers from expansion and eccentric endocardial displacement during systole. This eccentric displacement is both a cause for lower pumping efficacy but also an abnormal loading condition for the underlying tissue. During the subsequent several days, the affected tissue thins and expands in a process called infarct expansion. Infarct expansion represents the apparent widening of the infarcted region, not by extension to non-affected areas, but by the thinning of the wall and increase of the surface of the affected part of the chamber~\cite{weisman_myocardial_1987}. While cell necrosis occurs within the infarct zone, and the necrotic tissue reabsorption leads to wall thinning, some experiments have found that the infarct expansion and wall thinning occurs even when no significant necrotic tissue reabsorption has yet occurred~\cite{hochman_expansion_1982}. This indicates that the wall thinning in infarct expansion may not only be due to mass loss. On the other hand, cell slippage, leading to a rearrangement of structures comprising the myocardium, was measured to occur in the expanded myocardium in rat models~\cite{weisman_cellular_1988}. Infarct expansion is linked to collagen degradation in the infarcted area~\cite{ whittaker_role_1991,hwang_situ_2017}. Collagen degradation is mainly speculated to be an active process, in the sense that it arises through biochemical signalling, due to the activation of matrix metalloproteinases within the infarcted region~\cite{cleutjens_regulation_1995}. After the first week, and during the next couple of weeks, the local collagen production increases in the infarcted region. Initially, this is a compensatory process which stiffens the affected tissue and limits the eccentric systolic stretch. A few weeks after the infarct, a stiff collagen scar appears in the affected region and the entire affected chamber has dilated~\cite{sutton_left_2000}.
\par
The mathematical modelling of growth and remodelling of soft biological tissues mainly follows two approaches: the theory of evolving reference configurations, also called multiplicative or kinetic growth theory, and the theory of constrained mixtures~\cite{ambrosi_growth_2019}. The theory of evolving reference configurations~\cite{goriely_mathematics_2017} has taken from plasticity theory the decomposition of the full deformation gradient into an elastic deformation tensor and an inelastic part, with the inelastic part possibly having determinant different from one. In this way, mass loss or gain can be modelled by having the determinant of the inelastic deformation be less than, or alternately greater than, one. The theory of evolving reference configurations has been incorporated in various works on the growth and remodelling of the myocardium~\cite{saez_computational_2016,avazmohammadi_interactions_2019,goktepe_multiscale_2010,goktepe_generic_2010, sharifi_multiscale_2024,himpel_computational_2005}. The evolving reference configuration methodology fails to capture changes at the constitutive level, since the constitutive law remains the same. The gross mechanical properties of the body change only due to the shape reprogramming, the added or lost stiffness of the body due to the rescaling of the free energy by the determinant of the growth tensor, and possible geometrical frustrations induced by incompatible growth. Therefore, where the underlying structure of the material changes, the evolving reference configurations is not sufficient to capture remodelling. Moreover, the choice for the evolution equations remains phenomenological. For example, in~\cite{goktepe_generic_2010}, the authors use the same critical pressure criterion for two different growth laws, an isotropic one and a transversely isotropic one. From a modelling standpoint this is unsatisfactory since the choice of growth seems chosen a priori and not a consequence of the biological and mechanical factors at play. Moreover, as in~\cite{goktepe_multiscale_2010}, the bounds for the growth do not result through the achievement of some biomechanical equilibrium, but are a model parameter, indicating that the growth evolution does not balance its driving force, whether it be strain or stress. The multiplicative growth decomposition has also been applied to the modelling of acute myocardial infarction in~\cite{saez_computational_2016}. In this model, infarct expansion occurs through thinning by way of cellular death and radial expansion due to the collagen degradation. However, the collagen degradation is not modelled explicitly, as the constitutive model remains the same. Moreover, there is no material softening due to the collagen degradation; the softening that occurs is due to the explicit mass loss due to the determinant of the growth tensor being smaller than one. Indeed, if the determinant were equal to one, due to the inelastic deformation in the presumably stiffer direction in the fibres, the material would actually harden.
\par
Constrained mixture theory has also been a popular approach for the modelling of soft tissues~\cite{humphrey_constrained_2002,gebauer_homogenized_2023,cyron_homogenized_2016,gierig_computational_2021,gierig_arterial_2023}. Homogenized constrained mixture theory with respect to the classical multiplicative growth theory treats each constituent of the material separately. Similarly to the multiplicative growth theory, the remodelling of the tissue is captured by assuming a different inelastic pre-stretch to each constituent. Therefore, the limitations attributed to the kinematic growth theory concerning the phenomenological choices of the inelastic deformations and their evolutions pertain also to the homogenized constrained mixture theory. Homogenized constrained mixture models typically feature a large number of parameters which need to be calibrated while the biomechanical properties of each individual constituent are not necessarily known. It is more flexible for the treatment of materials composed of multiple constituents, which is often the case for soft tissues.
\par
We propose instead to model the remodelling by including changes in the stored energy function directly.  We believe that damage can be regarded as a type of remodelling as it affects the constitutive law governing the solid. It describes the failure of a material due to the accumulation of microscopic failure. Damage can be linked to the formation of microscopic voids and cracks and is an irreversible process~\cite{lemaitre_mechanics_1990}. According to the theory proposed by Griffith~\cite{Griffith1921}, fracture and damage in a material evolve according to a principle of minimal energy. Indeed, if a body's equilibrium configuration is such that it minimizes its total mechanical energy, then the equilibrium can be a damaged state, if such state can be reached by a continuous evolution decreasing the total mechanical energy. It provides a very clear and universal first principle according to which inelastic phenomena in a body occur, that is in contrast with the modelling approaches described above, whose evolution laws are mainly phenomenological in nature. Damage is an alternative to both the kinematic growth theory and the constrained mixture theory, as the reference unstressed configuration does not evolve in time. Its evolution is governed by the Griffith's total energy minimization principle. Moreover, damage can be coupled to plasticity models describing inelastic deformations~\cite{lemaitre_coupled_1985,marino_molecular-level_2019}. 
\par
In this work, we propose a model derived from variational principles for the onset and evolution of the damage to the collagen scaffold within the myocardium in the first few days following a cardiac infarction. We start from the variational formulation of the damage problem, where the deformation and damage fields are minimizers of a particular total mechanical energy. From the necessary conditions such minimizers must satisfy, we derive a model for irreversible damage evolution. We specialise this model for the case of the damage of the collagen in the cleavage planes within the myocardium. We apply our model to two test cases. The first test case is the indentation of a passive rectangular geometry, reproducing the systolic stretch post myocardial infarction, where we find that increased local load inducing an excessive shearing of the sheetlets produces damage. The second test case concerns the damage onset post-myocardial infarction. Damage to the collagen scaffold connecting the myocytes in the infarcted area would explain the relative cardiomyocyte slippage seen in infarct expansion. Indeed, with no additional assumptions other than the loss of contractility in the infarcted area, we find that damage evolves within the infarcted region. To the best of our knowledge, this is the first work to describe post-myocardial infarction remodelling in damage mechanics setting. 
\par
This work is structured as follows. In Section~\ref{sec:models_and_methods} we give the elements of our mechanical model of postinfarction remodelling. We introduce mathematical notations and conventions used throughout this work and review damage modelling in a variational setting. Starting from the variational problem, we derive evolution equations for the damage field, and specialise the damage problem in the case of the myocardium. We propose a numerical scheme for the numerical approximation of the damage-electromechanics problem. In Section~\ref{sec:test_cases}, we present our two test cases and simulate the evolution of the damage in an indented slab geometry and in a left ventricle following a myocardial infarction. Finally, in Section~\ref{sec:discussion}, we summarise the main contributions, as well as their limitations and future developments.
\section{Models and methods}
\label{sec:models_and_methods}
\subsection{General notions of continuum mechanics}
\label{sec:continuum_mechanics}
Let us consider an unloaded and unstressed homogeneous body $\mathcal B$ occupying the open bounded domain immersed in flat space $\Hat{\Omega}\subset\mathbb R^3$. We call $\Hat{\Omega}$ the reference unstressed configuration of the body. We assume that the body is continuously deformed in time into the configurations $\mathcal B_t$ occupying the bounded domains $\Omega_t$ by a suitably smooth family of injective deformation maps $\bm{\chi}_t$ such that
\begin{equation*}
    \begin{aligned}
        \bm{\chi}_t:\Hat{\Omega}\ni\bm{X}\mapsto\bm x\in\Omega_t, \quad\forall t\in[0,T].
    \end{aligned}
\end{equation*}
In particular, we assume that the deformation maps can be expressed as:
\begin{equation*}
    \bm x = \bm{\chi}_t(\bm X) = \bm X+\bm d(\bm X;t),
\end{equation*}
where $\bm d$ is the displacement of a given point in the reference configuration. We denote with $\bm F$ the deformation gradient, $\bm F = \nabla\bm{\chi}_t = \bm I+\nabla\bm d$, where $\bm I$ is the identity tensor, and with $\bm C$ we denote the right Cauchy-Green deformation tensor, $\bm C = \bm F^T\bm F$. We remark that the deformations $\bm{\chi}_t$ must be such that they preserve the orientation, i.e. $J=\mathrm{det}\bm F>0$. We call the set of all second-order tensors with positive determinant $\mathrm{Lin^+}$.
\par
We recall that a material is called hyperelastic, if there exists a free energy or strain energy density $\psi$ defined on the reference configuration $\Hat{\Omega}$,
\begin{equation*}
    \psi:\mathrm{Lin}^+\to[0,\infty),
\end{equation*}
such that for the first Piola-Kirchhoff stress tensor $\bm P= J\bm \sigma \bm F^{-T}$, where $\bm \sigma$ is the Cauchy stress tensor,
it holds, equivalently, that:
\begin{equation*}
    \bm P = \pdv{\psi}{\bm F} = 2\bm F\pdv{\psi}{\bm C},
\end{equation*}
where we use the following conventions for tensor derivatives:
\begin{equation*}
    \left(\pdv{\psi}{\bm F}\right)_{i,j} = \pdv{\psi}{ F_{i,j}},\quad \mbox{and}\quad \left(\nabla\cdot\bm P\right)_i = \pdv{P_{i,j}}{X_j}.
\end{equation*}
For a hyperelastic material, all closed strain cycles are non-dissipative.
\par
Let us assume that a body is loaded quasi-statically so that at each time $t_n\in[0,T]$ it is in equilibrium, and that the boundary of the reference configuration of the body is partitioned into disjoint open subsets $\Hat{\Gamma} _\mathrm{D}$ and $\Hat{\Gamma} _\mathrm{N}$, such that $\partial\Hat\Omega = \overline{\Hat{\Gamma}} _\mathrm{N} \cup \overline{\Hat{\Gamma}} _\mathrm{D}, \quad\Hat{\Gamma}_\mathrm{N} \cap \Hat{\Gamma} _\mathrm{D} = \emptyset, \quad \bm{\chi}_t(\Hat{\Gamma} _\mathrm{N}) = {\Gamma} _\mathrm{N}, \quad \bm{\chi}_t(\Hat{\Gamma} _\mathrm{D}) = {\Gamma} _\mathrm{D}$. At each quasi-static loading time $t_n$ the body is subject to a zero-displacement constraint on the Dirichlet boundary $\Hat{\Gamma} _\mathrm{D}$ and a nominal traction $\mathbf{g}_n$ in the reference configuration acting on the boundary $\Hat{\Gamma} _\mathrm{N}$. Assuming no body forces, the total mechanical energy of a deformed hyperelastic body at time $t_n$ subject to the displacement and traction constraints is:
\begin{equation}
\label{eq:total_energy}
    \mathcal{E}_n(\bm d) = \int_{\Hat{\Omega}}\psi(\bm F)\ \mathrm{d}\Hat{\Omega} - \int_{\Hat{\Gamma} _\mathrm{N}}\mathbf{g}_n\cdot\bm d\ \mathrm{d}\Hat{\Gamma}.
\end{equation}
We are interested in finding configurations $\Omega$ that minimize the mechanical energy of the body at each discrete time $t_n=n\Delta t$, for $n=0,1,\ldots N_t$, where $N_t$ is the number of time steps and $\Delta t>0$ the time step size. Due to the characterization of the deformation via displacements, this is equivalent to finding the displacements $\bm d_n\in V$ for which the energy functional~\eqref{eq:total_energy} is minimal:
\begin{equation}
    \label{eq:minima}
    \mathcal{E}_n(\bm d_n) = \min_{{\bm d}\in V} \mathcal{E}_n({\bm d}),
\end{equation}
for all $n=0,1,\ldots N_t$ where $V$ is the space of all admissible displacements,
\begin{equation*}
    V = [H^1_{ \Hat\Gamma_\mathrm{D}}(\Hat\Omega)]^3.
\end{equation*}
We recall that the displacement field $\bm d_n$ is a local minimum of the total mechanical energy functional~\eqref{eq:total_energy}, only if its Fréchet derivative is identically zero, i.e. if:
\begin{equation}
\label{eq:el_mechanics}
    \mathcal D_{\bm d}\mathcal E_n(\bm d_n)[\bm v] = \int_{\Hat\Omega}\pdv{\psi(\bm F_n)}{\bm F}:\nabla \bm v\ \mathrm{d}\Hat\Omega - \int_{\Hat\Gamma_\mathrm{N}}\mathbf{g}_n\cdot\bm v\ \mathrm{d}\Hat\Gamma = 0 \qquad \forall\bm v\in V.
\end{equation}
By testing Eq.~\eqref{eq:el_mechanics} with sufficiently regular admissible displacements, we obtain the Euler-Lagrange equations:
\begin{subequations}
    \label{eq:el_eqs}
    \begin{empheq}[left=\empheqlbrace]{align}
     & -\nabla\cdot\bm P_n = \bm 0 &&\mbox{in} \ \Hat\Omega, \label{eq:el_strong}
    \\
    & \bm d_n = \bm 0 &&\mbox{on} \ \Hat\Gamma_\mathrm{D}, \label{eq:bcs_dirichlet_strong}
    \\
    & \bm P_n\bm N = \mathbf{g}_n \hspace{0.75em} &&\mbox{on} \ \Hat\Gamma_\mathrm{N}, \label{eq:bcs_neumann_strong}
    \end{empheq}
\end{subequations}
for all $n=0,1,\ldots N_t$, where $\bm N$ is the outward normal unit vector to the body in the reference configuration. Eq.~\eqref{eq:el_strong} is the balance of linear momentum in local form written in the reference configuration. This means that finding equilibrium configurations for the body for a hyperelastic material reduces to finding minimizers to the total mechanical energy functional~\eqref{eq:total_energy}.
\subsection{Damage modelling}
\label{sec:damage_modelling}
Let us assume that the body $\mathcal{B}$ is prone to damage. We model damage as a type of material failure or dissipation under which the body does not behave elastically, i.e. under specific loading conditions the strain energy density function degrades. We introduce a damage variable
\begin{equation*}
    \alpha:\Hat\Omega\to [0,1],
\end{equation*}
where $\alpha(\bm X)=0$ if at the material point $\bm X$ no damage has occurred, whereas $\alpha(\bm X)=1$ if at the material point $\bm X$ the material has fully degraded. The constitutive law of the material changes according to the damage variable $\alpha$,
\begin{equation*}
    \partial_\alpha\psi_\mathrm{d}(\bm F,\alpha) \leq 0\quad \forall\bm F\in\mathrm{Lin}^+,
\end{equation*}
so that as damage increases, the stored energy decreases and the material softens.
For a class of materials including damage called gradient damage materials~\cite{pham_gradient_2011}, the total mechanical energy of the system at time $t_n$, subject to displacement and traction constraints, depends on the damage variable $\alpha$, and reads\footnote{For a symmetric positive-definite second-order tensor $\bm K$, the tensor $\sqrt{\bm K}$ is defined as the only symmetric positive-definite second-order tensor such that $\sqrt{\bm K}^T\sqrt{\bm K} = \bm K$.}:
\begin{equation}
    \label{eq:energy_gradient_damage}
    \mathcal{E}_{\ell,n}(\bm d,\alpha) = \int_{\Hat{\Omega}}\psi_\mathrm{d}(\bm F,\alpha)\ \mathrm{d}\Hat{\Omega} + \int_{\Hat{\Omega}}w_1\left[\psi_w(\alpha)+\ell^2|\sqrt{\bm K}\nabla\alpha|^2\right]\ \mathrm{d}\Hat{\Omega} - \int_{\Hat{\Gamma} _\mathrm{N}}\mathbf{g}_n\cdot\bm d\ \mathrm{d}\Hat{\Gamma},
\end{equation}
where $\bm K$ is a symmetric positive-definite second-order tensor, $w_1\geq 0$ is a specific fracture energy parameter, $\psi_w(\alpha)\geq 0$ is a dissipated energy function, such that
\begin{equation*}
    \psi_w:[0,1]\to[0,1], \quad \psi_w'(\alpha)\geq 0, \quad \psi_w(0) = 0,\ \mbox{and} \quad \psi_w(1) = 1,
\end{equation*}
and $\ell>0$ is a material length-scale parameter. We consider set the dissipated energy function to be:
\begin{equation}
\label{eq:dissipation_def}
    \psi_w(\alpha) = \alpha^2,
\end{equation}
in accordance with the linear elastic case~\cite{jean-jacques_marigo_gradient_2016}. Following~\cite{bourdin_numerical_2000}, the specific fracture energy parameter $w_1$ can be expressed as:
\begin{equation}
\label{eq:constant_fracture_energy}
    w_1 = \frac{G_c}{2\ell},
\end{equation}
in light of its connection with the phase-field modelling of fracture, where $G_c$ is a specific surface fracture energy.
In particular, we model the hyperelastic energy as:
\begin{equation}
\label{eq:psi_mod_def}
    \psi_\mathrm{d}(\bm F,\alpha) = \psi(\bm F) - g(\alpha)\psi_\mathrm{diss}(\bm F),
\end{equation}
where, following~\cite{miehe_phase_2010,borden_phase-field_2012,borden_phase-field_2016}, $g(\alpha)$ is the degradation function such that:
\begin{equation*}
    g:[0,1]\to[0,1], \quad g'(\alpha)\geq 0, \quad g(0) = 0,\ \mbox{and} \quad g(1) = 1,
\end{equation*}
which, following~\cite{jean-jacques_marigo_gradient_2016} for the linear elastic case, we assume to be:
\begin{equation}
\label{eq:degradation_def}
    g(\alpha) = 1 - (1-\alpha)^2,
\end{equation}
and $\psi_\mathrm{diss}(\bm F)$ is a dissipated energy,
\begin{equation*}
    \psi_\mathrm{diss}:\mathrm{Lin}^+\to[0,\infty).
\end{equation*}
In order for the energy $\psi_\mathrm{d}(\bm{F},\alpha)$ to be well-defined in Eq.~\eqref{eq:psi_mod_def}, it has to hold that
\begin{equation}
\label{eq:well_posed_dissipation_energy}
    \psi_\mathrm{diss}(\bm F)\leq\psi(\bm F), \quad \forall\bm F\in\mathrm{Lin}^+.
\end{equation}
We will assume that the damage is irreversible, that is, if $\alpha_{n-1}$ and $\alpha_n$ are the damage fields at times $t_{n-1}$ and $t_n$ respectively, for all $n=1,\ldots,N_t$, the following irreversibility constraint holds:
\begin{equation}
\label{eq:damage_irr_constraint}
    \alpha_n\geq\alpha_{n-1}\ \text{a.e. in }\hat{\Omega}.
\end{equation}
Let $W$ be the space of admissible damage fields:
\begin{equation*}
    W = H^1_{\Hat\Gamma_\alpha}(\Hat{\Omega}),
\end{equation*}
where $\Hat\Gamma_\alpha\subseteq\partial\Omega$ is the Dirichlet boundary for the damage $\alpha$, where we assume no damage occurs. Similarly to problem~\eqref{eq:minima}, at each time $t_n\in[0,T]$, we are looking for an admissible displacement $\bm d_n$ and a damage field $\alpha_n$ that minimize the total energy functional~\eqref{eq:energy_gradient_damage} under the irreversibility constraint~\eqref{eq:damage_irr_constraint}:
\begin{equation*}
    \mathcal{E}_{\ell,n}(\bm d_n,\alpha_n) = \min_{\substack{{\bm d}\in V,\ {\alpha}\in W\\ {\alpha}\in [0,1] \\ \alpha\geq \alpha_{n-1}\ \text{a.e. in }\hat{\Omega}}}\mathcal{E}_{\ell,n}({\bm d},{\alpha}), \quad \forall n=0,\ldots,n.
\end{equation*}
We will subsequently drop the boundedness constraint $\alpha\in[0,1]$ as it is automatically satisfied for a suitable extension of the degradation function $g$ to $\mathbb R$, as explained in~\ref{appendix:boundedness}.
\par
Under assumptions of sufficient smoothness of the strain energy density, the minimizers  $\bm d_n,\alpha_n$ of the total mechanical energy~\eqref{eq:energy_gradient_damage} satisfy the following Euler-Lagrange equation for the displacement:
\begin{equation}
    \label{eq:el_damage_1}
    \mathcal{D}_{\bm d}\mathcal{E}_{\ell,n}(\bm d_n,\alpha_n)[\delta \bm d] = \int_{\Hat\Omega}\pdv{\psi_\mathrm{d}(\bm F_n,\alpha_n)}{\bm F}:\nabla \delta\bm d\ \mathrm{d}\Hat\Omega - \int_{\Hat\Gamma_\mathrm{N}}\mathbf{g}_n\cdot\delta\bm d\ \mathrm{d}\Hat\Gamma = 0, \quad \forall\delta\bm d\in V,\quad \forall n=1,\ldots N_t,
\end{equation}
which is similar to Eq.~\eqref{eq:el_mechanics}, where the damage variable $\alpha_n$ acts as a softening material parameter.
\par
On the other hand, given the damage at the previous time $\alpha_{n-1}$, the minimizing damage field $\alpha_n$ must satisfy the following necessary Karush-Kuhn-Tucker conditions for all $n=0,\ldots,N_t$:
\begin{subequations}
    \label{eq:kkt}
    \begin{empheq}[left=\empheqlbrace]{align}
    & \mathcal{D}_\alpha\mathcal{E}_{\ell,n}(\bm d_n, \alpha_n)[\delta\alpha] = \int_{\Hat{\Omega}}\mu_n\delta\alpha\  \mathrm{d}\Hat\Omega, \quad \forall \delta\alpha\in H_{\Hat\Gamma_\alpha}^1(\Hat\Omega), \label{eq:kkt_stability}\\
    & \alpha_{n}\geq\alpha_{n-1},\quad \mu_n\geq 0,\quad \mu_n(\alpha_{n}-\alpha_{n-1})=0,\quad \mbox{a.e. in }\Hat{\Omega},\label{eq:kkt_energy balance}
    \end{empheq}
\end{subequations}
for some $\mu_n\in L^2(\Hat{\Omega})$.
\par
A strategy for finding local minima of the total mechanical energy, and thus solutions to the damage problem $\bm d_n,\alpha_n$, would be solving Eq.~\eqref{eq:el_damage_1} and Eq.s~\eqref{eq:kkt}. However, Eq.s~\eqref{eq:kkt} are not straightforwardly solvable using standard Galerkin finite element methods. This drives us to propose an alternative to Eq.s~\eqref{eq:kkt} for a governing equation for the damage variable $\alpha_n$, one which may be solved using standard Galerkin finite element methods~\cite{quarteroni_numerical_2017}, while preserving the irreversibility of the damage evolution. We proceed to derive such a governing equation, starting from~\eqref{eq:kkt}.
\par
First of all, we remark that, given $\bm d_n$, if the unique $\alpha_n$ such that
\begin{equation}
\label{eq:pseudo_el_damage}
\begin{aligned}
    \mathcal{D}_\alpha\mathcal{E}_{\ell,n}(\bm d_n, \alpha_n)[\delta\alpha]=&\int_{\Hat{\Omega}} - g'(\alpha_n)\psi_\mathrm{diss}(\bm{F}_n)v\ \mathrm{d}\Hat\Omega\ +\ \int_{\Hat{\Omega}} w_1{\psi_w'(\alpha_n)}v\ \mathrm{d}\Hat\Omega\ \\&+\ \int_{\Hat{\Omega}} 2w_1\ell^2{\bm K}\nabla\alpha_n\cdot\nabla v\ \mathrm{d}\Hat\Omega = 0, \qquad\forall v \in H^1_{\Hat\Gamma_\alpha}(\Hat{\Omega}),
\end{aligned}
\end{equation}
is such that the irreversibility constraint~\eqref{eq:damage_irr_constraint} is satisfied, then conditions~\eqref{eq:kkt} are satisfied with $\mu=0$. The uniqueness of $\alpha_n$ for Eq.~\eqref{eq:pseudo_el_damage} is guaranteed by the Lax-Milgram lemma~\cite{evans_partial_2022}, since thanks to the definition of $g(\alpha)$ and $w(\alpha)$, Eq.~\eqref{eq:pseudo_el_damage} is linear in $\alpha$. Let $\alpha_{n-1},\alpha_{n}$ be damage fields, satisfying Eq.~\eqref{eq:pseudo_el_damage} at times $t_{n-1}$ and ${t_{n}}$, respectively, and let $\delta\alpha = \alpha_{n}-\alpha_{n-1}$. Then, recalling the definitions of $\psi_w(\alpha)$ and $g(\alpha)$ in~\eqref{eq:dissipation_def} and~\eqref{eq:degradation_def}:
\begin{equation*}
\begin{aligned}
    & \int_{\Hat{\Omega}} 2(\alpha_{n}\psi_\mathrm{diss}(\bm{F}_{n})-\alpha_{n-1}\psi_\mathrm{diss}(\bm{F}_{n-1}))v\  \mathrm{d}\Hat\Omega+\int_{\Hat{\Omega}} 2w_1\delta\alpha v\  \mathrm{d}\Hat\Omega\ \\ &+\ \int_{\Hat{\Omega}} 2w_1\ell^2{\bm K}\nabla\delta\alpha\cdot\nabla v\ \mathrm{d}\Hat\Omega = \int_{\Hat{\Omega}} 2(\psi_\mathrm{diss}(\bm{F}_{n}) -\psi_\mathrm{diss}(\bm{F}_{n-1})) v\ \mathrm{d}\Hat\Omega, \quad\forall v \in H^1_{\Hat\Gamma_\alpha}(\Hat{\Omega}),\quad\forall n=1,\ldots,N_t,
\end{aligned}
\end{equation*}
which yields:
\begin{equation*}
\begin{aligned}
    & \int_{\Hat{\Omega}} \psi_\mathrm{diss}(\bm{F}_{n})\delta\alpha v\  \mathrm{d}\Hat\Omega+\int_{\Hat{\Omega}} w_1\delta\alpha v\ \mathrm{d}\Hat\Omega+\int_{\Hat{\Omega}} w_1\ell^2{\bm K}\nabla\delta\alpha\cdot\nabla v\ \mathrm{d}\Hat\Omega\ \\ & = \int_{\Hat{\Omega}}(1-\alpha_{n-1}) (\psi_\mathrm{diss}(\bm{F}_{n}) -\psi_\mathrm{diss}(\bm{F}_{n-1})) v\ \mathrm{d}\Hat\Omega, \quad\forall v \in H^1_{\Hat\Gamma_\alpha}(\Hat{\Omega}),\quad\forall n=1,\ldots,N_t.
\end{aligned}
\end{equation*}
Since it holds that $0\leq\alpha_{n-1}\leq 1$, if the dissipation energy were monotonically increasing in time, 
\begin{equation*}
    \psi_\mathrm{diss}(\bm{F}_{n}) \geq \psi_\mathrm{diss}(\bm{F}_{n-1}),\quad \forall n=1,\ldots,N_t, \quad \mbox{a.e. in }\Hat{\Omega},
\end{equation*}
then this would imply by the weak maximum principle~\cite{evans_partial_2022} that $\delta\alpha\geq 0$, ensuring the satisfaction of the irreversibility constraint~\eqref{eq:damage_irr_constraint}. This motivates the use of a history variable, as in~\cite{miehe_phase_2010}, as a driver for damage evolution. In particular, to ensure the irreversibility constraint to hold, we replace $\psi_\mathrm{diss}(\bm F_n)$ appearing in Eq.~\eqref{eq:pseudo_el_damage} by the history variable $\xi_n$, defined as:
\begin{equation}
\label{eq:def_truncated_history_variable}
    \xi_{n}(\psi_\mathrm{diss}) = \max_{k\in \{0,\dots,n\}}\psi_\mathrm{diss}(\bm{F}_{k}),
\end{equation}
with $\xi_0=0$. The history variable $\xi_n$ defined in Eq.~\eqref{eq:def_truncated_history_variable} is increasing in time, and therefore yields an irreversible damage field evolution.
\par
In the following, as the equation governing the damage evolution, we will be using:
\begin{equation}
\begin{aligned}
    \label{eq:damage_evolution}
    \int_{\Hat{\Omega}} - g'(\alpha_n)\xi_n\delta\alpha\ \mathrm{d}\Hat\Omega\ +\ \int_{\Hat{\Omega}} w_1{\psi_w'(\alpha_n)}\delta\alpha\ \mathrm{d}\Hat\Omega\ &+\ \int_{\Hat{\Omega}} 2w_1\ell^2\bm K\nabla\alpha_n\cdot\nabla \delta\alpha\ \mathrm{d}\Hat\Omega = 0, \\ &\forall \delta\alpha \in H^1_{\Hat\Gamma_\alpha}(\Hat{\Omega}),\quad \forall n=1,\ldots,N_t,
\end{aligned}
\end{equation}
with $\xi_n$ defined in~\eqref{eq:def_truncated_history_variable}.

\subsection{Myocardial constitutive modelling}
\label{sec:myocardial_structure_energetics}
Due to the characteristic structure of the myocardium and its distinct mechanical and electrical properties along different directions, we introduce a local reference system in accordance with the principal material directions.\par
In particular, for any $\bm X\in\Hat{\Omega}$, we define an orthonormal basis $\{\hat{\bm f},\hat{\bm s},\hat{\bm n}\}$ of the tangent space at $\bm X$, which reflects the principal material properties, representing the principal local fibre direction $\hat{\bm f}$, the sheet direction $\hat{\bm s}$, and the sheet-normal direction $\hat{\bm n}$. In particular, $\hat{\bm f}\cross\hat{\bm s}$ generates the surface element of the sheetlet. We say the basis reflects the principal properties insofar cardiomyocytes have a natural dispersion around a particular direction.
\par
The basis $\{\hat{\bm f},\hat{\bm s},\hat{\bm n}\}$ rotates along $\hat{\Omega}$. Furthermore, upon an arbitrary deformation, the basis is pushed forward into a no longer orthonormal triplet $\{\bm f,\bm s, \bm n\}$, where:
\begin{equation*}
    \bm f = \frac{\bm F\bm \Hat{\bm f}}{\|\bm F\bm \Hat{\bm f}\|},\quad \bm s = \frac{\bm F\bm \Hat{\bm s}}{\|\bm F\bm \Hat{\bm s}\|}, \quad \bm n = \frac{\bm F\bm \Hat{\bm n}}{\|\bm F\bm \Hat{\bm n}\|}
\end{equation*}
where $\bm f$ is the principal direction of the fibres in the current configuration $\Omega$.
\par
Due to the symmetries in the microstructure, the material can be assumed to be orthotropic, with distinct material properties in the three directions $\{\hat{\bm f},\hat{\bm s},\hat{\bm n}\}$. A widely accepted model for the passive myocardium, which replicates available experimental data, and whose expression highlights the dependence of the material response on the principal material directions $\{\hat{\bm f},\hat{\bm s},\hat{\bm n}\}$ is the exponential incompressible Holzapfel-Ogden model~\cite{holzapfel_constitutive_2009}. In particular, for the undamaged part of the elastic energy in Eq.~\eqref{eq:psi_mod_def}, we use a modified quasi-incompressible Holzapfel-Ogden model~\cite{holzapfel_constitutive_2009, nolan_robust_2014}, where we write explicitly the exponential stiffening in the $\hat{\bm n}$ direction, as:
\begin{equation}
    \label{eq:ho_mod_qinc}
    \begin{aligned}
        \psi(\bm F) &= \frac{a}{2b}\exp[b(\Bar{I}_1-3)] + \sum_{i=f,s,n}\frac{a_i}{2b_i}\left\{\exp[b_i\langle I_{4i}-1\rangle_+^2]-1 \right\} \\ &+ \frac{a_{fs}}{2b_{fs}}\left[\exp(b_{fs}I_{8fs}^2)-1\right] + \frac{c_\mathrm{bulk}}{2}\left(J-1\right)\mathrm{log}(J),
    \end{aligned}
\end{equation}
where $I_1$ and $I_{4f},I_{4s},I_{8fs}$ are respectively the first isochoric isotropic invariant, and the fourth and eighth anisotropic invariants of the right Cauchy-Green tensor $\bm C$, defined as:
\begin{equation*}
    \begin{aligned}
        \Bar{I}_1 &= \mathrm{tr}( J^{-2/3}{\bm C}),\\
        I_{4i} &= \Hat{\bm i}\cdot\bm C\Hat{\bm i},\quad && i\in\{f,s,n\} \\
        I_{8ij} &= \Hat{\bm i}\cdot\bm C\Hat{\bm j},\quad && i,j\in\{f,s,n\},\quad i\neq j,
    \end{aligned}
\end{equation*}
and $\langle x\rangle_+$ denotes the positive part of $x$. We set the parameters of the strain-energy density~\eqref{eq:ho_mod_qinc} as to replicate experimental shear data from~\cite{dokos_shear_2002}, by considering analytical solutions to the homogeneous simple shear problem as in~\cite{holzapfel_constitutive_2009}.
\par
Regarding the dissipated energy term $\psi_\mathrm{diss}$, we assume damage to occur to the collagen fibre family connecting the sheetlets, which is principally oriented in the $\hat{\bm n}$ direction. In line with other material models of fibre-reinforced soft biological tissues~\cite{holzapfel_new_2000}, we assume that the dissipated part $\psi_\mathrm{diss}$ of the elastic energy~\eqref{eq:psi_mod_def} is:
\begin{equation}
    \label{eq:psi_diss_def}
    \psi_\mathrm{diss}(\bm F) = \frac{a_n}{2b_n}\left\{\exp[b_n\langle I_{4n}-1\rangle_+^2]-1\right\},
\end{equation}
which trivially satisfies~\eqref{eq:well_posed_dissipation_energy}. Due to the dissipated term~\eqref{eq:psi_diss_def} not being explicitly modelled in the Holzapfel-Ogden strain energy density~\cite{holzapfel_constitutive_2009}, the resulting model's stresses may no longer be monotonous with respect to the strains. Therefore, the modification of the Holzapfel-Ogden model~\eqref{eq:ho_mod_qinc} by adding explicit dependence on the invariant $I_{4n}$ is necessary. Further details are given in~\ref{appendix:monotonicity_of_the_stress}.
\par
We model the damage that occurs in to the collagen strands in the cleavage planes separating the sheetlets made of parallel myocytes. We assume the damage propagates more easily along than across cleavage planes. We recall the sheetlets and thus the cleavage planes to be orthogonal to $\hat{\bm n}$. We include this transverse isotropy by defining the damage propagation tensor $\bm K$ appearing in~\eqref{eq:damage_evolution} as:
\begin{equation}
\label{eq:k_spectral_form}
    \bm K = k\Hat{\bm f}\otimes\Hat{\bm f} + k\Hat{\bm s}\otimes\Hat{\bm s} + \Hat{\bm n}\otimes\Hat{\bm n},
\end{equation}
with $k>1$ being a material parameter which we assume to be constant in space.
\par
Because of the electrically induced cardiac contraction, the myocardium is considered an active material. We consider the active stress approach~\cite{ambrosi_active_2012} to model the microscopic stress generated by the myofilaments, by which the total Piola-Kirchhoff stress tensor $\bm P$ may be separated into a passive contribution $\bm P^\mathrm{pass}$, and an active part $\bm P^{\mathrm{act}}$ as:
\begin{equation*}
    \bm P = \bm P^\mathrm{pass} + \bm P^\mathrm{act},
\end{equation*}
which in the case of hyperelastic materials specializes to:
\begin{equation*}
    \bm P(\bm d, \alpha, f_\mathrm{isch}) = \pdv{\psi_\mathrm{d}}{\bm F}\left(\bm F, \alpha\right) + \bm P^\mathrm{act}(\bm F,f_\mathrm{isch}).
\end{equation*}
We assume the cardiac contraction occurs as consequence of a space and time-dependent active tension generated along the fibres in the reference configuration:
\begin{equation}
\label{eq:ta_ischemia}
    \bm P^\mathrm{act}(\bm F_n,f_{\mathrm{isch},n}) = f_\mathrm{isch}(\bm X,t_n)\Hat{T}_\mathrm{a}(t_n)\frac{\bm F_n\Hat{\bm f}\otimes\Hat{\bm f}}{\|\bm F_n\Hat{\bm f}\|},
\end{equation}
where $\Hat{T}_\mathrm{a}(t)$ is the active tension generated at a generic healthy point $\bm X\in\Hat{\Omega}$, and with
\begin{equation*}
    f_\mathrm{isch}:\Hat{\Omega}\times[0,T]\to[0,1]
\end{equation*}
being a known function modelling the ischaemic effects, with $\partial_tf_\mathrm{isch}(\bm X,t)\leq 0$, such that $f_\mathrm{isch}(\bm X,t)= 0$ if at point $\bm X$ there is no oxygenated blood supply, $f_\mathrm{isch}(\bm X,t)= 1$ if at $\bm X$ the tissue is fully healthy, and $f_\mathrm{isch}(\bm X,t)\in (0,1)$ if $\bm X$ is at the infarct border zone.

\par
The heartbeat is an extremely dynamic event, therefore a quasi-static approximation for the loading conditions during the whole heartbeat is not appropriate. We can extend the model from Section~\ref{sec:damage_modelling} to the dynamic time-continuous case by assuming the damage to be instantaneous with respect to the displacement and coupling directly the time-continuous representation of Eq.~\eqref{eq:damage_evolution} with the elastodynamics equation. In particular, the more general time-continuous damage mechanics problem in a reference domain $\hat{\Omega}$ consists of finding the pair $(\bm d(t),\alpha(t))$, for all $t\in(0,T]$, such that:
\begin{subequations}
\label{eq:damage_mechanics_unspec_bcs}
    \begin{empheq}[left=\empheqlbrace]{align}
    & \rho\pdv[2]{\bm{d}}{t} -\nabla\cdot\left(\bm P^\mathrm{pass}\left(\bm F, \alpha\right) + \bm P^\mathrm{act}  (\bm F,f_{\mathrm{isch}})\right) = \bm{0}, && \ \mbox{in}\ \Hat{\Omega}\cross(0,T],
    \label{eq:mechanics_strong}
    \\
    &- g'(\alpha)\xi(t)+w_1{\psi_w'(\alpha)}- \nabla\cdot\left(2w_1\ell^2\bm K\nabla\alpha\right) = 0, && \ \mbox{in}\ \Hat{\Omega}\cross(0,T],
    \label{eq:damage_strong}
    \\
    &\bm d = \bm 0, && \ \mbox{on}\ \Hat{\Gamma}_\mathrm{D}\cross(0,T],
    \label{eq:mech_bc_1}
    \\
    & \bm P\bm N = \mathbf{g}(t), && \ \mbox{on}\ \Hat{\Gamma}_\mathrm{N}\cross(0,T],
    \label{eq:mech_bc_2}
    \\
    & \alpha = 0, && \ \mbox{on}\ \Hat{\Gamma}_\alpha\cross(0,T],
    \label{eq:damage_bc_1}
    \\
    &\bm K\nabla\alpha\cdot\bm N = 0, && \ \mbox{on}\ \partial\Hat{\Omega}\setminus\Hat{\Gamma}_\alpha\cross(0,T],
    \label{eq:damage_bc_2}
    \\
    & \left(\bm{d},\dot{\bm{d}},\alpha\right)(0)=\left(\bm{d}_0,\bm{0},0\right), &&\ \mbox{in} \ \Hat{\Omega},
    \label{eq:ic}
    \end{empheq}
\end{subequations}
with $\bm F = \bm I+\nabla\bm d$ and $\xi(t)$ defined as:
\begin{equation*}
    \label{eq:}
    \xi(t) = \max_{t\in (0,T]}\psi_\mathrm{diss}(\bm{F}(t)),
\end{equation*}
 with $\xi(0)=0$. Eq.~\eqref{eq:mechanics_strong} is the momentum conservation equation and taking into consideration also the active mechanical behaviour of the myocardium, Eq.~\eqref{eq:damage_strong} is the damage evolution equation, Eqs.~\eqref{eq:mech_bc_1} and~\eqref{eq:mech_bc_1} are the Dirichlet and Neumann boundary conditions for~\eqref{eq:mechanics_strong}, while Eqs.~\eqref{eq:damage_bc_1} and~\eqref{eq:damage_bc_2} are the homogeneous Dirichlet and Neumann boundary conditions for~\eqref{eq:damage_strong}, and~\eqref{eq:ic} is the initial condition. We will specify the boundary conditions further in Sections~\ref{sec:slab} and~\ref{sec:acute_infarct}.
\subsection{A note on material homogeneity and residual stresses}
\label{sec:residual_stresses_strains}
The left ventricle is residually stressed~\cite{omens_residual_1990}. This is a common occurrence in biological tissues, and can be seen using the opening-angle experiment. In this experiment, circumferential slices of the ventricles, arteries, and other specimens are radially cut and subsequently opened~\cite{omens_residual_1990,holzapfel_layer-specific_2007}, exposing the presence residual circumferential and axial residual stress in the intact samples. Moreover, in~\cite{holzapfel_layer-specific_2007}, the authors affirm that, at least in the aorta, the residual stresses are of a three-dimensional nature an cannot be reduced to a single parameter. The presence of residual stresses in excised biological tissues indicates that neither the imaging configuration nor the excised configuration represent the reference configuration, which is by definition stress-free.
\par
Some of these issues have been tackled separately in previous works by the solid mechanics community. In~\cite{barnafi_reconstructing_2024}, the authors recover the reference configuration from an imaging configuration under recorded under a given load by formulating an inverse problem in the displacement. Due to the formulation of the inverse problem, the obtained reference configurations are virtual only due to the global incompatibility constraint. In~\cite{genet_heterogeneous_2015}, the authors introduce a method for the inclusion of heterogeneous isotropic growth-induced prestrain in the left ventricle starting from an ungrown unloaded configuration. Still, from a generic imaging configuration, there is no unique way of introducing a prestrain measure on the heart.
\par
Residual stresses in the heart are non-negligible and may be of different origin. Moreover the strain energy distribution over the single heartbeat depends on the definition of the reference configuration and therefore on the modelled residual stresses. In our analysis, we include the effect of residual stresses and material inhomogeneities in the left ventricle by modelling an explicitly space-dependent specific fracture energy parameter $w_1(\bm X)$. In particular, given a space-dependent reference dissipation energy $\psi_{\mathrm{diss},\mathrm{ref}}(\bm X)$, we define the space-dependent specific fracture energy parameter as:
\begin{equation}
\label{eq:specific_fracture_energy_space_dep}
    w_1(\bm X) = \frac{G_c}{2\ell}\exp\left(\frac{\|\psi_\mathrm{diss,ref}(\bm X, \cdot)\|_{L^\infty([0,T])}}{\|\psi_\mathrm{diss,ref}(\cdot,\cdot)\|_{L^\infty(\Hat\Omega\times[0,T])}}\right) + \|\psi_\mathrm{diss,ref}(\bm X, \cdot)\|_{L^\infty([0,T])}.
\end{equation}
In particular, we define the reference dissipation energy $\psi_\mathrm{diss,ref}(\bm X, t)$ as:
\begin{equation*}
    \psi_\mathrm{diss,ref}(\bm X, t) =  \psi_\mathrm{diss}({\bm F}(\hat{\bm d}(t)),t),
\end{equation*}
where $\hat{\bm d}$ is the displacement solution to Eq.~\eqref{eq:mechanics_strong} given $f_\mathrm{isch}\equiv 1$ and $\alpha\equiv 0$, i.e. for a fully healthy left ventricle. The specific fracture energy $w_1(\bm X)$ given by~\eqref{eq:specific_fracture_energy_space_dep} means that fracture occurs only if the reference energy is exceeded by a factor ranging in $(\frac{G_c}{2\ell},e\frac{G_c}{2\ell})$, and that to create damage in points which normally carry higher energies more energy needs to be invested.
\subsection{Numerical approximation}
\label{sec:num_approximation}
We present the numerical approximation for Eqs.~\eqref{eq:damage_mechanics_unspec_bcs}. Since Eq.~\eqref{eq:mechanics_strong} with fixed $\alpha_n$ is the classic elastostatics problem, and Eq.~\eqref{eq:damage_strong} with $\xi_n$ fixed is a linear reaction-diffusion equation, we apply a partitioned scheme to solve Eq.s~\eqref{eq:mechanics_strong} and~\eqref{eq:damage_strong} separately. In particular, we apply a staggered, loosely coupled scheme. In the case of the coupling with a cardiac electromechanics model for the recovery of the active tension $\hat{T}_\mathrm{a}$ in Eq.~\eqref{eq:ta_ischemia}, the damage-mechanics problem is solved following this scheme, while the rest of the problem is solved following~\cite{stella_fast_2022,regazzoni_cardiac_2022}.
\par
Regarding the spatial approximation of~\eqref{eq:damage_mechanics_unspec_bcs}, we use continuous nodal finite elements. We introduce the mesh $\mathcal T_h$ of the computational domain $\Hat\Omega$ made of tetrahedral elements with average maximal diameter $h$. We denote with $\mathbb P^1(K)$ the set of linear polynomials over the mesh element $K\in\mathcal T_h$. We introduce the following finite-dimensional subspaces $V_h$ and $W_h$:
\begin{equation*}
    \begin{aligned}
        &V_h = \left\{\bm v_h\in [H^1(\Hat{\Omega})]^3:\ \bm v_{h\vert_K}\in [\mathbb P^1(K)]^3,\ \forall K\in\mathcal T_h\right\}, \\
        &W_h = \left\{w_h \in H^1(\Hat{\Omega}):\  w_{h\vert_K}\in \mathbb P^1(K),\ \forall K\in\mathcal T_h\right\},
    \end{aligned}
\end{equation*}
for which we introduce their nodal bases, $\left\{\bm \varphi_i\right\}_{i=1}^{\mathrm{dim}V_h}$ and $\left\{ \varphi_i\right\}_{i=1}^{\mathrm{dim}W_h}$ respectively. Regarding the time discretization, we use finite differences with a step size of $\Delta t$.
\par
Given the initial damage field $\alpha_{h,0}=0$, the staggered Galerkin finite element formulation of problem~\eqref{eq:damage_mechanics_unspec_bcs} reads as follows. At time $t_n=n\Delta t$, $n=1,\ldots N_t$, given the active stress $T_\mathrm{a}(t_n)$, the boundary loads $\mathbf{g}_n$, and $\bm d_{n-1}, \bm d_{n-2},\alpha_{h,n-1},\xi_{h,n-1} $, find $\bm d_{h,n}\in V_h$ and $\alpha_{h,n}\in W_h$, such that $\bm d_{h,n} = \bm 0$ on $\Hat\Gamma_\mathrm{D}$, $\alpha_{h,n} = 0$ on $\Hat\Gamma_\alpha$ and:
\begin{equation}
\begin{aligned}
    \int_{\Hat\Omega}\rho\frac{\bm d_{h,n}-2\bm d_{h,n-1}+\bm d_{h,n-2}}{\Delta t^2}\cdot \bm \varphi_i\ \mathrm{d}\Hat\Omega &+ \int_{\Hat\Omega}\left(\bm P^\mathrm{pass}(\bm F_{h,n},\alpha_{h,n-1})+\bm P^\mathrm{act}(\bm F_{h,n},T_\mathrm{a}(t_n))\right):\nabla \bm \varphi_i\ \mathrm{d}\Hat\Omega \\ &- \int_{\Hat\Gamma_\mathrm{N}}\mathbf{g}_n\cdot\bm \varphi_i\ \mathrm{d}\Hat\Gamma = 0, \quad \forall i\in\left\{1,\dots,\mathrm{dim}V_h\right\},
\end{aligned}
\end{equation}
and
\begin{equation}
\begin{aligned}
    &\int_{\Hat\Omega}\left( -g'(\alpha_{h,n})\xi_{h,n}(\psi_{\mathrm{diss},h,n})+w_1\psi_w'(\alpha_{h,n}) \right)\varphi_i\ \mathrm{d}\Hat\Omega \\ &+\int_{\Hat\Omega}2w_1\ell^2\bm K\nabla\alpha_{h,n}\cdot\nabla\varphi_i\ \mathrm{d}\Hat\Omega = 0, \quad \forall i\in\left\{1,\dots,\mathrm{dim}W_h\right\}.
\end{aligned}
\end{equation}
The quantities $\bm F_{h,n}$ and $\psi_{\mathrm{diss},h,n}\in W_h$ are such that:
\begin{equation*}
    \bm F_{h,n} = \bm I+\nabla\bm d_{h,n},
\end{equation*}
and
\begin{equation*}
\int_{\Hat\Omega}\psi_{\mathrm{diss},h,n}\varphi_i\mathrm{d}\Hat\Omega = \int_{\Hat\Omega}\psi_{\mathrm{diss}}(\bm F_{h,n})\varphi_i\mathrm{d}\Hat\Omega + \int_{\Hat\Omega}h_K^2\nabla\psi_{\mathrm{diss}}(\bm F_{h,n})\cdot\nabla\varphi_i\ \mathrm{d}\Hat\Omega,\quad \forall i\in\left\{1,\dots\mathrm{dim}W_h\right\}.
\end{equation*}

\section{Test cases}
\label{sec:test_cases}
We apply our model in two test cases. First, in Section~\ref{sec:slab}, we study the damage in a passive rectangular slab. Then, in Section~\ref{sec:acute_infarct}, we model the onset of damage following an acute myocardial infarction in a realistic left ventricular geometry.
\begin{figure}[t]
     \centering
     \subfloat[][]{\def\svgwidth{5.5cm}
\begingroup%
  \makeatletter%
  \providecommand\color[2][]{%
    \errmessage{(Inkscape) Color is used for the text in Inkscape, but the package 'color.sty' is not loaded}%
    \renewcommand\color[2][]{}%
  }%
  \providecommand\transparent[1]{%
    \errmessage{(Inkscape) Transparency is used (non-zero) for the text in Inkscape, but the package 'transparent.sty' is not loaded}%
    \renewcommand\transparent[1]{}%
  }%
  \providecommand\rotatebox[2]{#2}%
  \newcommand*\fsize{\dimexpr\f@size pt\relax}%
  \newcommand*\lineheight[1]{\fontsize{\fsize}{#1\fsize}\selectfont}%
  \ifx\svgwidth\undefined%
    \setlength{\unitlength}{536.70885942bp}%
    \ifx\svgscale\undefined%
      \relax%
    \else%
      \setlength{\unitlength}{\unitlength * \real{\svgscale}}%
    \fi%
  \else%
    \setlength{\unitlength}{\svgwidth}%
  \fi%
  \global\let\svgwidth\undefined%
  \global\let\svgscale\undefined%
  \makeatother%
  \begin{picture}(1,0.97652733)%
    \lineheight{1}%
    \setlength\tabcolsep{0pt}%
    \put(0,0){\includegraphics[width=\unitlength,page=1]{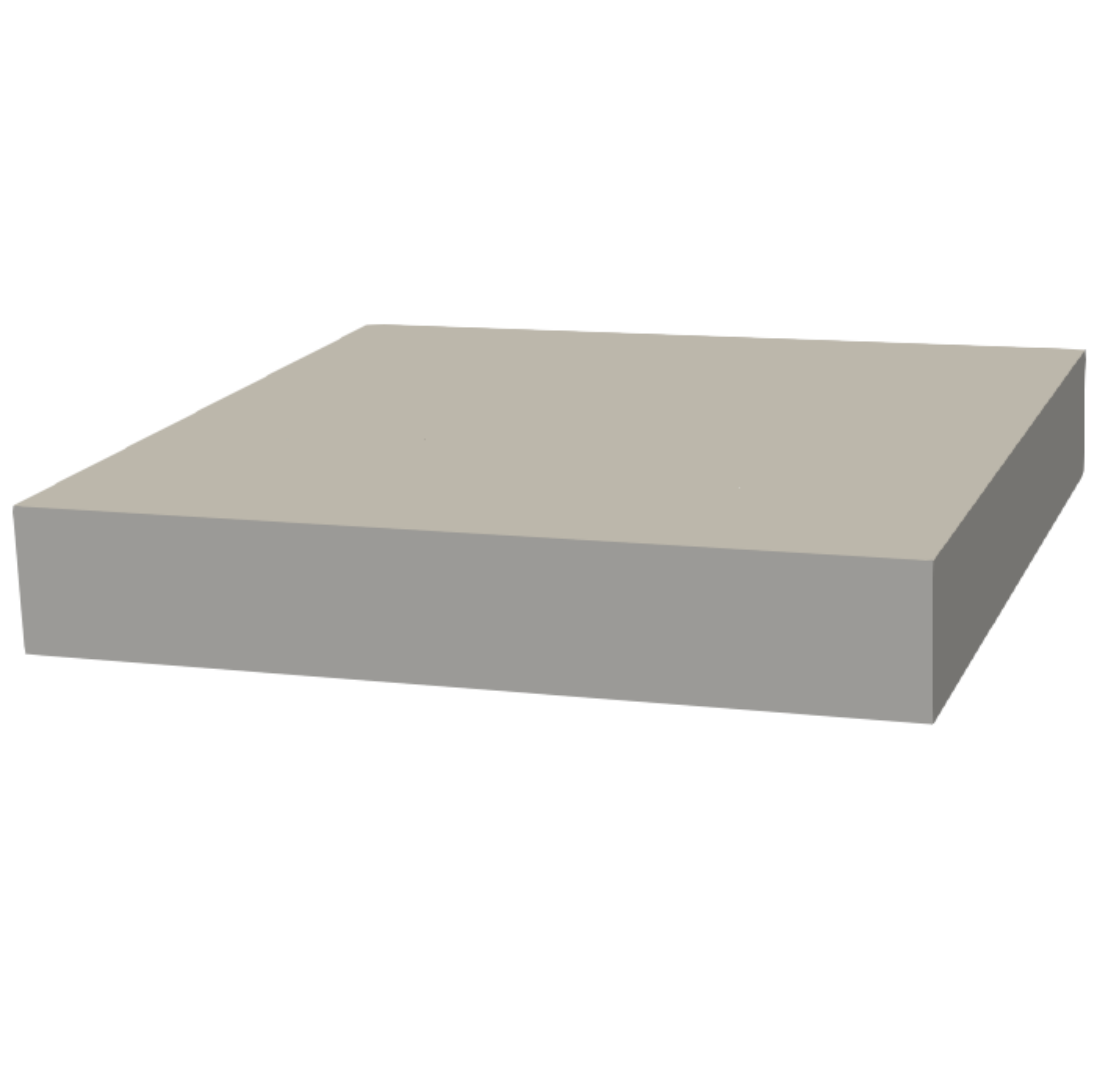}}%
    \put(0.92565159,0.30674459){\makebox(0,0)[t]{\lineheight{1.25}\smash{\begin{tabular}[t]{c}$\Hat\Omega$\end{tabular}}}}%
    \put(0.51035849,0.30827483){\makebox(0,0)[t]{\lineheight{1.25}\smash{\begin{tabular}[t]{c}$\Hat\Gamma^\mathrm{epi}$\end{tabular}}}}%
    \put(0.93303902,0.50728418){\makebox(0,0)[t]{\lineheight{1.25}\smash{\begin{tabular}[t]{c}$\Hat\Gamma^\mathrm{side}$\end{tabular}}}}%
    \put(0.51573093,0.59052989){\makebox(0,0)[t]{\lineheight{1.25}\smash{\begin{tabular}[t]{c}$\Hat\Gamma^\mathrm{endo}$\end{tabular}}}}%
    \put(0,0){\includegraphics[width=\unitlength,page=2]{slab_domain_boundaries.pdf}}%
  \end{picture}%
\endgroup%
\label{fig:slab_domain_boundaries}}
     \qquad\qquad
     \subfloat[][]{\def\svgwidth{5.5cm}
\begingroup%
  \makeatletter%
  \providecommand\color[2][]{%
    \errmessage{(Inkscape) Color is used for the text in Inkscape, but the package 'color.sty' is not loaded}%
    \renewcommand\color[2][]{}%
  }%
  \providecommand\transparent[1]{%
    \errmessage{(Inkscape) Transparency is used (non-zero) for the text in Inkscape, but the package 'transparent.sty' is not loaded}%
    \renewcommand\transparent[1]{}%
  }%
  \providecommand\rotatebox[2]{#2}%
  \newcommand*\fsize{\dimexpr\f@size pt\relax}%
  \newcommand*\lineheight[1]{\fontsize{\fsize}{#1\fsize}\selectfont}%
  \ifx\svgwidth\undefined%
    \setlength{\unitlength}{467.25bp}%
    \ifx\svgscale\undefined%
      \relax%
    \else%
      \setlength{\unitlength}{\unitlength * \real{\svgscale}}%
    \fi%
  \else%
    \setlength{\unitlength}{\svgwidth}%
  \fi%
  \global\let\svgwidth\undefined%
  \global\let\svgscale\undefined%
  \makeatother%
  \begin{picture}(1,1.13804173)%
    \lineheight{1}%
    \setlength\tabcolsep{0pt}%
    \put(0,0){\includegraphics[width=\unitlength,page=1]{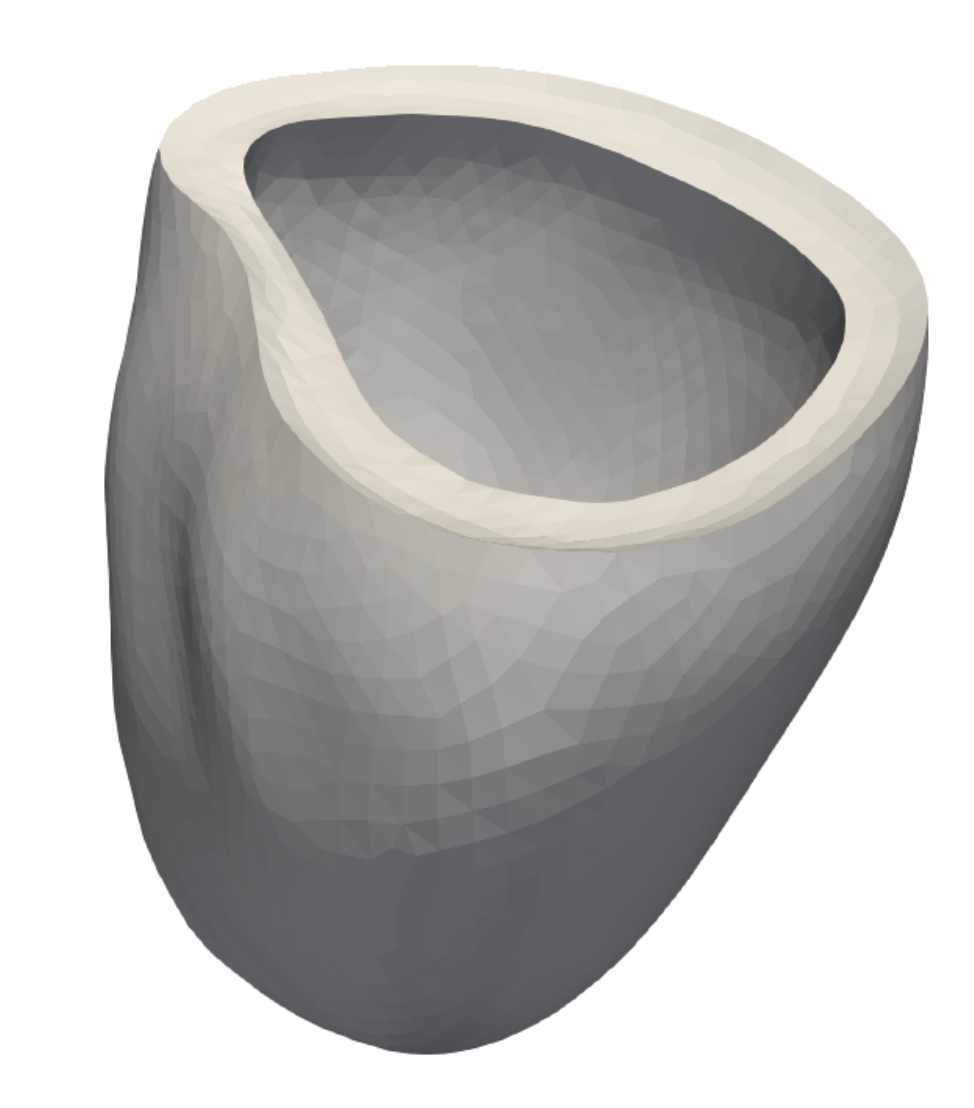}}%
    \put(0.6548165,0.07267475){\makebox(0,0)[t]{\lineheight{1.25}\smash{\begin{tabular}[t]{c}$\Hat\Omega$\end{tabular}}}}%
    \put(0.53723757,0.43236526){\makebox(0,0)[t]{\lineheight{1.25}\smash{\begin{tabular}[t]{c}$\Hat\Gamma^\mathrm{epi}$\end{tabular}}}}%
    \put(0.25835267,1.03405765){\makebox(0,0)[t]{\lineheight{1.25}\smash{\begin{tabular}[t]{c}$\Hat\Gamma^\mathrm{base}$\end{tabular}}}}%
    \put(0.57525864,0.80359541){\makebox(0,0)[t]{\lineheight{1.25}\smash{\begin{tabular}[t]{c}$\Hat\Gamma^\mathrm{endo}$\end{tabular}}}}%
  \end{picture}%
\endgroup%
\label{fig:lv_domain_boundaries}}
     \caption{{Reference domains. (a): Slab reference configuration $\Hat\Omega$, along with its boundaries $\Hat\Gamma^\mathrm{endo}$, $\Hat\Gamma^\mathrm{epi}$, $\Hat\Gamma^\mathrm{side}$. (b): Representation of the reference configuration $\Hat\Omega$ of the left ventricle, given by the Zygote Solid 3D left ventricle \cite{zygote_media_group_inc_zygote_2014}, along with its boundaries $\Hat\Gamma^\mathrm{endo}$, $\Hat\Gamma^\mathrm{epi}$, $\Hat\Gamma^\mathrm{base}$.}}
\end{figure}

\subsection{Test 1: Indentation of a passive slab geometry}
\label{sec:slab}
We study the damage induced by an indenter-like load applied on a passive rectangular slab geometry, depicted in Figure~\ref{fig:slab_domain_boundaries}. The local increase of passive loading mimics the increased load in the infarcted area found following a myocardial infarction. The slab has a thickness of $1~\si{\centi\metre}$ along the $z$ direction and a width of $6~\si{\centi\metre}$ along the $x$ and $y$ directions. We assume the load to be applied slowly, so we have that at each time the solid is in equilibrium and so we may disregard the intertial term in Eq.~\eqref{eq:mechanics_strong}. The boundary conditions for the slab are the following:
\begin{subequations}
    \label{eq:slab-system}
    \begin{empheq}[left=\empheqlbrace]{align}
    & \bm P^\mathrm{pass}\left(\bm F, \alpha\right)\bm{N} = -p_\mathrm{ind}(\bm X,t)\bm{N}, && \mbox{on} \ \Hat{\Gamma}^\mathrm{endo}\cross(0,T], \label{eq:bc_endo_slab}
    \\
    & \bm P^\mathrm{pass}\left(\bm F, \alpha\right)\bm{N} + \mathbf{K}^\mathrm{epi}\bm{d} = \bm{0}, && \mbox{on} \ \Hat{\Gamma}^\mathrm{epi}\cross(0,T],
    \label{eq:bc_epi_slab}
    \\
    & \bm d = \bm 0, && \mbox{on} \ \Hat{\Gamma}^\mathrm{side}\cross(0,T],
    \label{eq:bc_side_slab}
    \\
    & \alpha = 0, && \mbox{on} \ \Hat{\Gamma}^\mathrm{side}\cross(0,T],
    \label{eq:bc_side_damage_slab}
    \\
    & \bm K\nabla{\alpha}\cdot \bm N =  0, && \mbox{on} \ \Hat{\Gamma}^\mathrm{endo}\cup\Hat{\Gamma}^\mathrm{epi}\cross(0,T].
    \label{eq:bc_epi_damage_slab}
    \end{empheq}
\end{subequations}
In particular, to replicate load conditions occurring in the myocardium, we apply homogeneous Dirichlet displacement boundary conditions on the lateral faces $\Hat{\Gamma}^\mathrm{side}$ of the slab~\eqref{eq:bc_side_slab}, Robin-like boundary conditions~\eqref{eq:bc_epi_slab} on $\Hat\Gamma^\mathrm{epi}$, mimicking the action of the pericardial sac~\cite{pfaller_importance_2019}, and an indenter-like pressure~\eqref{eq:bc_endo_slab} at the endocardial surface $\Hat\Gamma^\mathrm{endo}$. The indenter-like pressure $p_\mathrm{ind}$ is such that it reproduces an indenter being pushed into the myocardial tissue, and thus takes the form of:
\begin{equation}
    \label{eq:pressure_indenter}
    p_\mathrm{ind}(\bm X,t)=\begin{dcases}f_\mathrm{ind}(\bm X)\ p_\mathrm{ind,max}, & \mbox {if } t > t_\mathrm{fin}, \\ f_\mathrm{ind}(\bm X)\ p_\mathrm{ind,max}\frac{t-t_\mathrm{init}}{t_\mathrm{fin}-t_\mathrm{init}}, & \mbox {if } t_\mathrm{init}\le t \le t_\mathrm{fin}, \\ 0, & \mbox {if } t < t_\mathrm{init},\end{dcases}
\end{equation}
where $p_\mathrm{ind,max}$ is a maximal indenter pressure and where the spatial dependence on the pressure $p_\mathrm{ind}(\bm X,t)$ is through $f_\mathrm{ind}(\bm X)$ which takes the form of:
\begin{equation}
    \label{eq:pressure_indenter_space}
    f_\mathrm{ind}(\bm X)=\begin{dcases}1, & \mbox {if } d(\bm X,\bm X_\mathrm{ind};\bm N) < r_\mathrm{ind,int}, \\ \frac{r_\mathrm{ind,ext}-d(\bm X,\bm X_\mathrm{ind};\bm N)}{r_\mathrm{ind,ext}-r_\mathrm{ind,int}}, & \mbox {if } r_\mathrm{ind,ext} \geq d(\bm X,\bm X_\mathrm{ind};\bm N) \geq r_\mathrm{ind,int}, \\ 0, & \mbox {if } d(\bm X,\bm X_\mathrm{ind};\bm N) > r_\mathrm{ind,ext},\end{dcases}
\end{equation}
where $d(\bm X,\bm X_\mathrm{ind};\bm N)$ is the distance between the point $\bm X$ and the axis parallel with the surface normal $\bm N$ passing through $\bm X_\mathrm{ind}$, $ r_\mathrm{ind,int}$ is an internal radius and $r_\mathrm{ind,ext}$ is an external radius. We assume the slab to be behaving only passively, hence the balance of linear momentum~\eqref{eq:mechanics_strong} contains only passive stresses. We include no material inhomogeneities, therefore the specific fracture energy is constant in space and takes the value of Eq.~\eqref{eq:constant_fracture_energy}.
Model and simulation parameters are reported in Tables~\ref{table:strain_energy_densities_params} and~\ref{table:damage_params}.
\begin{table}[t]
    \centering 
    \begin{tabular}{ p{5em} p{2em} p{3em} p{3em} p{3em}}
    Parameter & Unit & $\psi_\mathrm{HO}$ & $\psi$ & $\psi_\mathrm{diss}$ \\ 
    \hline
    $a$ & $\si{\pascal}$ & 59 & 54 & 0 \\
    $b$ & $-$ & 8.023 & 2.223 & 0 \\ 
    $a_f$ & $\si{\pascal}$ & 18472 & 21072 & 0 \\ 
    $b_f$ & $-$ & 16.026 & 15.026 & 0 \\ 
    $a_s$ & $\si{\pascal}$ & 2421 & 5642 & 0 \\ 
    $b_s$ & $-$ & 11.12 & 12.62 & 0 \\
    $a_n$ & $\si{\pascal}$ & 0 & 2821 & 2821 \\
    $b_n$ & $-$ & 0 & 12.62 & 12.62 \\
    $a_{fs}$ & $\si{\pascal}$ & 216 & 432 & 0 \\
    $b_{fs}$ & $-$ & 11.436 & 11.436 & 0 \\
    $c_\mathrm{bulk}$ & $\si{\pascal}$ & 0 & 50000 & 0 \\
    \end{tabular}
    \\[10pt]
    \caption{Strain energy density parameters.}
    \label{table:strain_energy_densities_params}
\end{table}

\begin{table}[t]
    \centering 
    \begin{tabular}{ p{5em} p{5em} p{3em} p{10em}}
    Parameter & Unit & Test 1 & Test 2  \\ 
    \hline
    $G_c$ & $\si{\pascal\metre}$ & $43$ & $43$ \\
    $\ell$ & $\si{\metre}$ & $0.06$ & $0.06$ \\
    $k$ & $-$ & $3$ & $3$ \\
    $\rho$ & $\si[per-mode=reciprocal]{\kilogram\per\cubic\meter}$  & \num[exponent-product=\ensuremath{\cdot}]{1000} \\
    $K_\perp$ & $\si[per-mode=reciprocal]{\pascal\per\meter}$ & 200000 & 200000 \\
    $K_\parallel$ & $\si[per-mode=reciprocal]{\pascal\per\meter}$ & 20000 & 20000\\
    $C_\perp$ & $\si[per-mode=reciprocal]{\pascal\second\per\meter}$ & 20000 & - \\
    $C_\parallel$ & $\si[per-mode=reciprocal]{\pascal\second\per\meter}$ & 2000 & - \\
    $\bm X_\mathrm{pz}$ & $\si{\metre}$ & $-$ & $(0.0693174,1.332,0.03)$  \\ 
    $r_\mathrm{pz}$ & $\si{\metre}$ & $-$ & $0.032$  \\ 
    $d_\mathrm{pz}$ & $\si{\metre}$ & $-$ & $0.01$  \\ 
    $t_\mathrm{init}$ & $\si{\second}$ & $0$ & $0.8$ \\ 
    $t_\mathrm{fin}$ & $\si{\second}$ & $12$ & $1.6$  \\
    $\bm X_\mathrm{ind}$ & $\si{\metre}$ & $-$ & $(0,0,0)$  \\  
    $p_\mathrm{ind,max}$ & $\si{\kilo\pascal}$ & $320$ & $-$  \\
    $r_\mathrm{ind,int}$ & $\si{\metre}$ & $0.005$ & $-$  \\ 
    $r_\mathrm{ind,ext}$ & $\si{\metre}$ & $0.007$ & $-$  \\ 
    $T$ & $\si{\second}$ & $12$ & $3.2$ \\
    $h$ & $\si{\metre}$ & $0.003$ & $0.003$ \\
    $\Delta t$ & $\si{\second}$ & $0.002$ & $0.001$ \\ 
    \end{tabular}
    \\[10pt]
    \caption{Model and simulation parameters for the slab geometry test case in Section~\ref{sec:slab} and for the acute infarction of the left ventricle in Section~\ref{sec:acute_infarct}.}
    \label{table:damage_params}
\end{table}

\begin{figure}[t]
     \centering
     \subfloat[][]{\includegraphics[width = 0.28\textwidth]{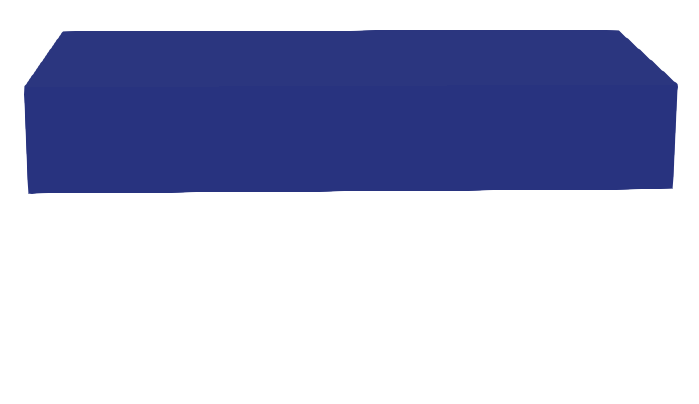}\label{fig:slab_0}}\hfill
     \subfloat[][]{\includegraphics[width = 0.28\textwidth]{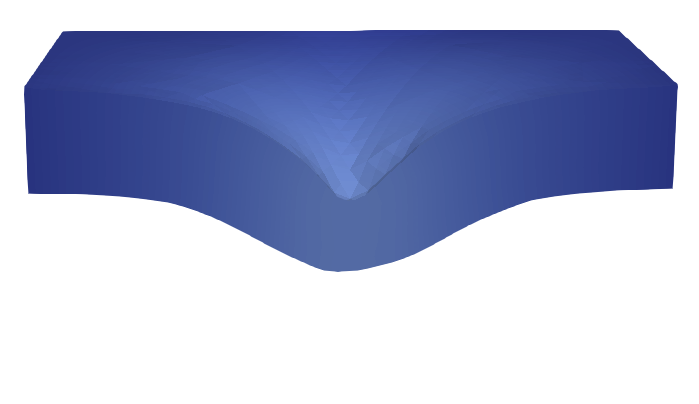}\label{fig:slab_1}}\hfill
     \subfloat[][]{\includegraphics[width = 0.28\textwidth]{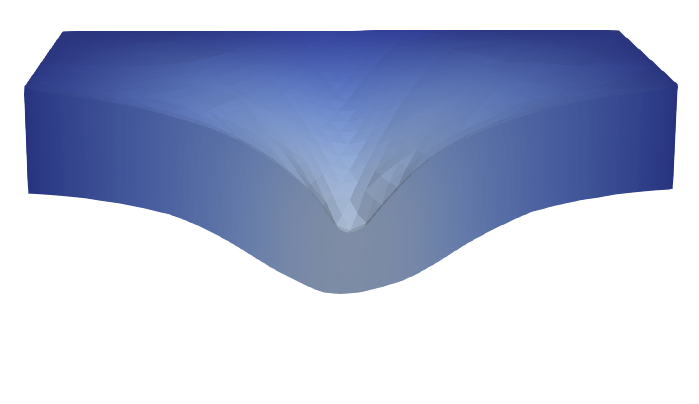}\label{fig:slab_2}}\\
     \subfloat[][]{\includegraphics[width = 0.28\textwidth]{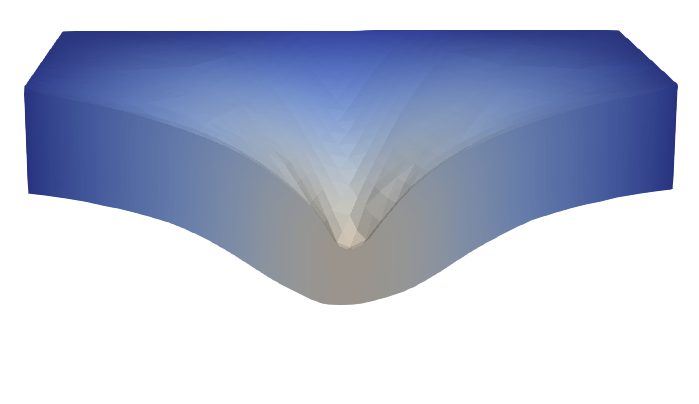}\label{fig:slab_3}}\hfill
     \subfloat[][]{\includegraphics[width = 0.28\textwidth]{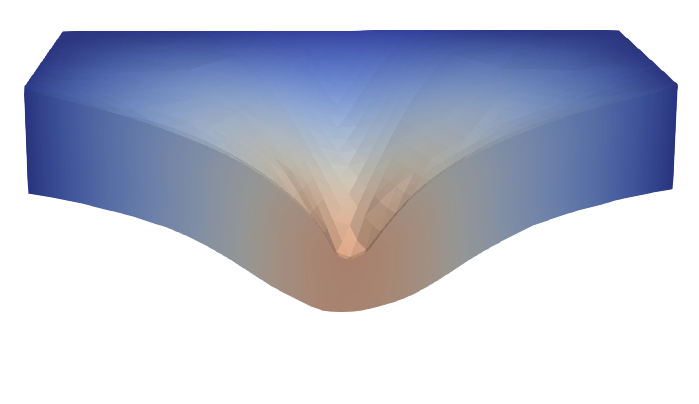}\label{fig:slab_4}}\hfill
     \subfloat[][]{\includegraphics[width = 0.28\textwidth]{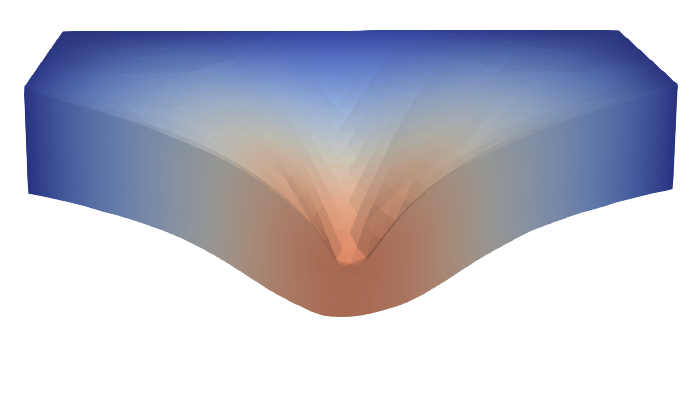}\label{fig:slab_5}}\\
     \subfloat[][]{\includegraphics[width = 0.28\textwidth]{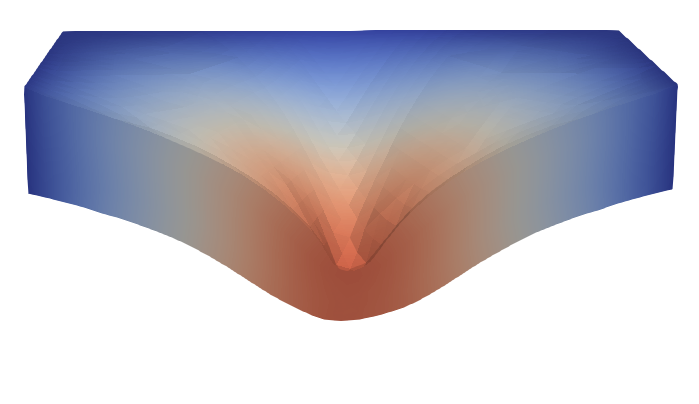}\label{fig:slab_6}}\hfill
     \subfloat[][]{\includegraphics[width = 0.28\textwidth]{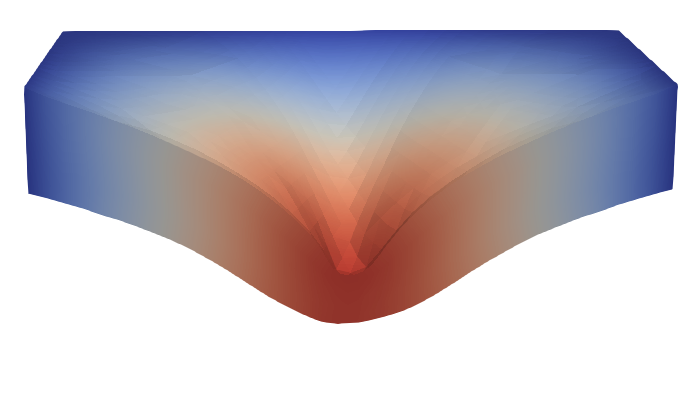}\label{fig:slab_7}}\hfill
     \subfloat[][]{\includegraphics[width = 0.28\textwidth]{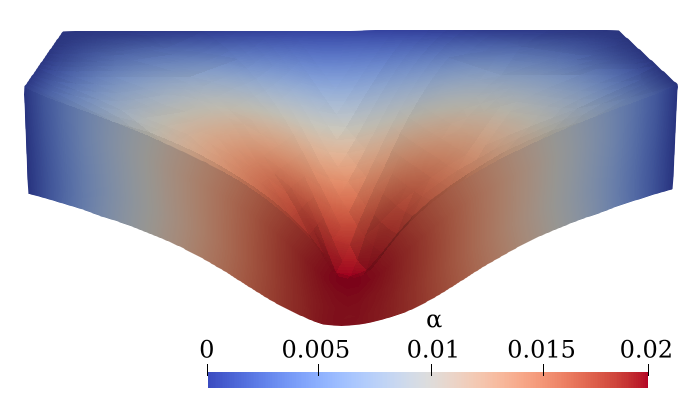}\label{fig:slab_8}}
     \caption{Test 1: Damage field $\alpha$ evolution on the cross-section perpendicular to the $y$ axis of the deformed geometry. (a) $t=0~\si{\second}$. (b) $t=1.52~\si{\second}$. (c) $t=3~\si{\second}$. (d) $t=4.52~\si{\second}$. (e) $t=6~\si{\second}$. (f) $t=7.52~\si{\second}$. (g) $t=9~\si{\second}$. (h) $t=10.52~\si{\second}$. (i) $t=12~\si{\second}$.}
     \label{fig:deformation_slab}
\end{figure}

\begin{figure}[t]
     \centering
     \subfloat[][]{\includegraphics[width = 0.28\textwidth]{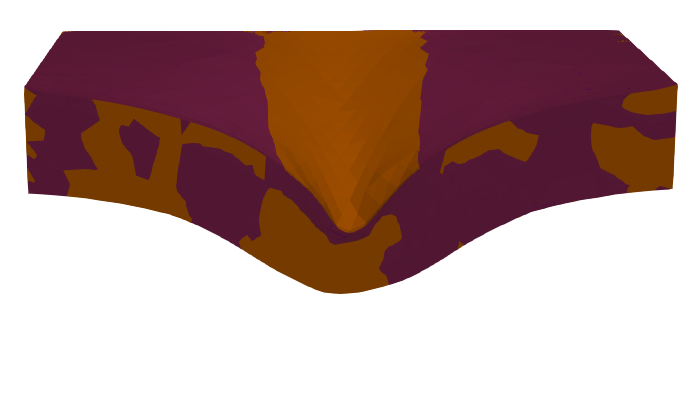}\label{fig:overlap1}}
     \subfloat[][]{\includegraphics[width = 0.28\textwidth]{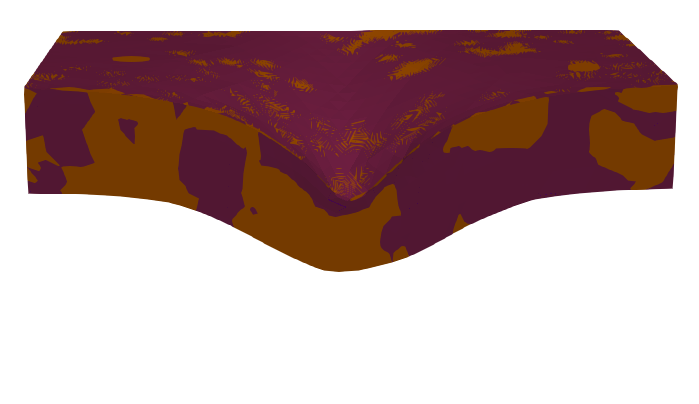}\label{fig:overlap2}}
     \caption{Test 1: Deformation of the cross-section perpendicular to the $y$ axis during loading and unloading. For the overlapping geometries the applied load is equal. The lighter colour corresponds to the configuration during loading, while the darker colour corresponds to the configuration during unloading. (a) $t=3~\si{\second}$ and $t=21~\si{\second}$. (b) $t=1.52~\si{\second}$ and $t=22.48~\si{\second}$.}
     \label{fig:deformation_slab_def}
\end{figure}

The fibres within the slab are oriented such that $\hat{\bm f}$ is aligned with the $y$ axis, $\hat{\bm s}$ is aligned with the $z$ axis and $\hat{\bm n}$ is aligned with the $x$ axis. The slab is loaded slowly with an increasing pressure according to~\eqref{eq:pressure_indenter}. The cross-section along the $z$ axis of the deformed slab along with the damage field is presented in Figure~\ref{fig:deformation_slab}.
From Figure~\ref{fig:deformation_slab}, we see that damage occurs instantly as the load is applied. Moreover the damage appears to be monotonously increasing, and occurring principally where the pressure is applied. Greater damage gradients are visible in the $\hat{\bm f}$ and  $\hat{\bm s}$ directions, as induced by the spectral form of the tensor $\bm K$ defined in Eq.~\eqref{eq:k_spectral_form}. In Figure~\ref{fig:deformation_slab_def}, we report the deformed damaged slab at equal indenter pressure during the loading and during a following unloading of the same duration.
From Figure~\ref{fig:deformation_slab_def}, we can see the effect of the damage on the deformation; at equal loads, during the unloading process, the tissue exhibits larger deformation. This is more evident at higher loads, reported in Figure~\ref{fig:overlap1}, where in the damaged region the deformation in the loading direction of the tissue during unloading seems to be overall greater. At smaller loadings, this is less evident, as exhibited by Figure~\ref{fig:overlap2}. This effect would be more pronounced if the damage field were greater.
\par
\begin{figure}[t]
     \centering
     \subfloat[][]{\includegraphics[width = 0.28\textwidth]{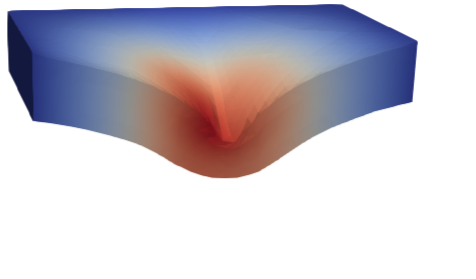}\label{fig:slabb_04}}\quad
     \subfloat[][]{\includegraphics[width = 0.28\textwidth]{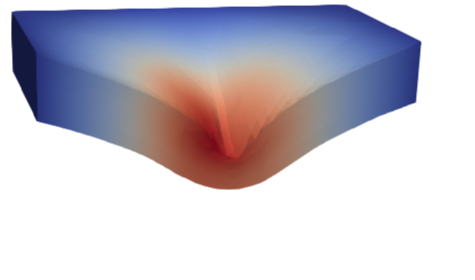}\label{fig:slabb_08}}\quad
     \subfloat[][]{\includegraphics[width = 0.28\textwidth]{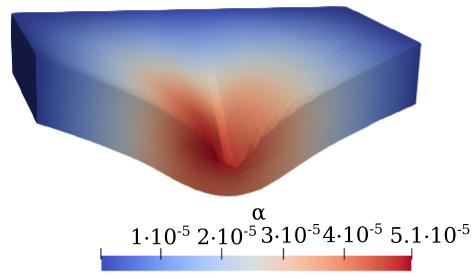}\label{fig:slabb_12}}\\
     \caption{Test 1: Damage field $\alpha$ on the cross-section perpendicular to the $y$ axis of the deformed geometry. (a) $t=4~\si{\second}$. (b) $t=8~\si{\second}$. (c) $t=12~\si{\second}$.}
     \label{fig:deformation_slabb}
\end{figure}

The orientation of the fibres within the geometry is crucial for the material response and the prediction of damage onset and evolution. When the fibres within the slab are oriented such that $\hat{\bm f}$ is aligned with the $x$ axis, $\hat{\bm s}$ is aligned with the $y$ axis and $\hat{\bm n}$ is aligned with the $z$ axis, the damage produced by the indenter-like load is significantly smaller, as evidenced by Figure~\ref{fig:deformation_slabb}. This is because, in the first case, the sheetlets spanned by $\hat{\bm f}$ and $\hat{\bm s}$ are oriented in such a way, that when the load is applied, the sheetlets shear one respect to the other, therefore causing an increase in the dissipation energy $\psi_\mathrm{diss}$. The damage field evolution in the second case happens early on, and the amount of damage does not visibly increase from Figure~\ref{fig:slabb_04} to Figure~\ref{fig:slabb_12} as the load increases. Furthermore, in Figure~\ref{fig:deformation_slabb} the damage field is not homogeneous across the thickness of the sample, but is nearer to the point of application of the load, while in Figure~\ref{fig:deformation_slab} the damage appears to be homogenous across the thickness of the sample.
\afterpage{\clearpage}
\subsection{Test 2: Acute infarction in a left ventricle}
\label{sec:acute_infarct}
\begin{figure}[t]
     \centering
     \subfloat[][]{\includegraphics[width = 0.2\textwidth]{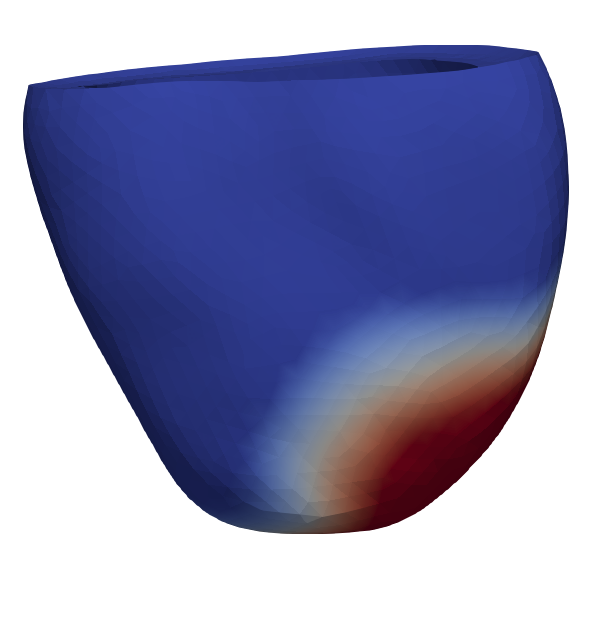}\label{fig:infarct_full}}\qquad
     \subfloat[][]{\includegraphics[width = 0.2\textwidth]{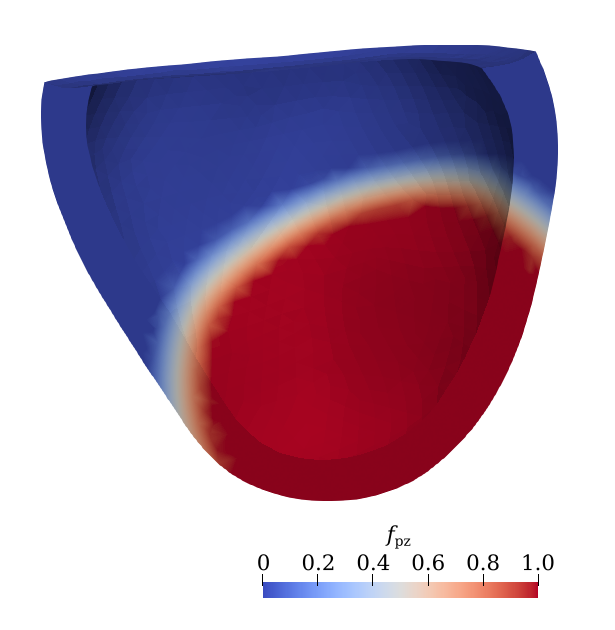}\label{fig:infarct_slice}}\\
     \caption{Infarcted portion of the myocardium on $\hat{\Omega}$. Values  of $f_\mathrm{pz}=1$ represent complete loss of contractility, whereas values of $f_\mathrm{pz}=0$ represent no contractility loss. (a) Complete view of the left ventricle. (b) Sliced view of the left ventricle.}
     \label{fig:infarcted_portion}
\end{figure}

\begin{figure}[t]
     \centering
     \subfloat[][]{\includegraphics[width = 0.2\textwidth]{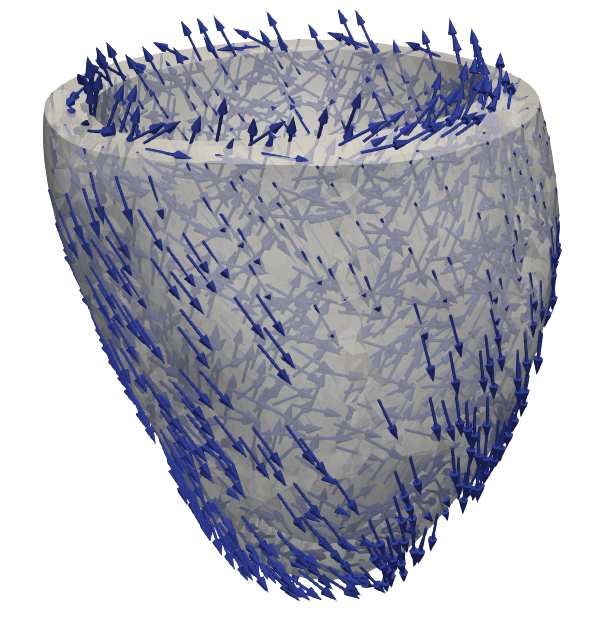}\label{fig:fibres_f}}\quad
     \subfloat[][]{\includegraphics[width = 0.2\textwidth]{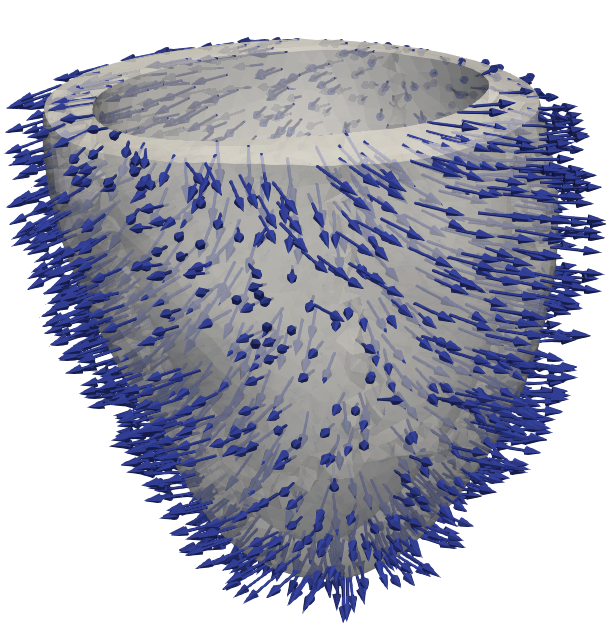}\label{fig:fibres_s}}\quad
     \subfloat[][]{\includegraphics[width = 0.2\textwidth]{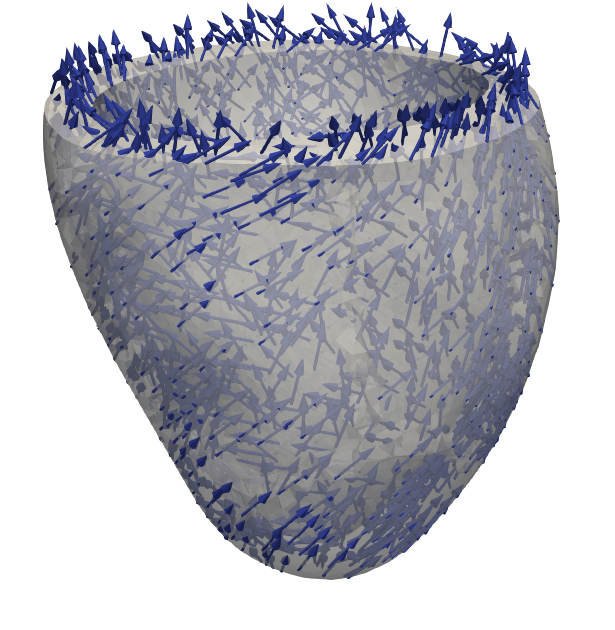}\label{fig:fibres_n}}\\
     \caption{Fibre orientation in the initial configuration with $\bm d = \bm d_0$. (a) ${\bm f}_0$ direction. (b) ${\bm s}_0$ direction. (c) ${\bm n}_0$ direction.}
     \label{fig:fibres_directions}
\end{figure}

We model the onset of damage to the collagen fibres in the cleavage planes within the myocardium, as consequence of an acute cardiac infarction. We consider the left ventricle as depicted in Figure~\ref{fig:lv_domain_boundaries}. The boundaries of the domain correspond to the endocardial surface $\hat{\Gamma}^\mathrm{endo}$, the epicardial surface $\hat{\Gamma}^\mathrm{epi}$, and the artificial boundary $\hat{\Gamma}^\mathrm{base}$, and are reported in Figure~\ref{fig:lv_domain_boundaries}. The boundary conditions for the damage-electromechanics model are:
\begin{subequations}
    \label{eq:damage_em_bcs}
    \begin{empheq}[left=\empheqlbrace]{align}
    & \bm P(\bm d, \alpha, f_\mathrm{isch})\bm{N} = -p_\mathrm{LV}(t)J\bm{F}^{-T}\bm{N}, && \mbox{on} \ \Hat{\Gamma}^\mathrm{endo}\cross(0,T], \label{eq:bc_endo}
    \\
    & \bm P(\bm d, \alpha, f_\mathrm{isch})\bm{N} + \mathbf{K}^\mathrm{epi}\bm{d} + \mathbf{C}^\mathrm{epi}\pdv{\bm{d}}{t} = \bm{0}, && \mbox{on} \ \Hat{\Gamma}^\mathrm{epi}\cross(0,T],
    \label{eq:bc_epi}
    \\
    & \bm P(\bm d, \alpha, f_\mathrm{isch})\bm{N} = p_\mathrm{LV}(t)||J\bm{F}^{-T}\bm{N}|| \frac{\int_{\Hat{\Gamma}^\mathrm{endo}}J\bm{F}^{-T}\bm{N} \mathrm{d}A }{\int_{\Hat\Gamma^\mathrm{base}}||J\bm{F}^{-T}\bm{N}||\mathrm{d}A}, && \mbox{on} \ \Hat{\Gamma}^\mathrm{base}\cross(0,T],
    \label{eq:bc_base}
    \\
    & \alpha = 0, && \mbox{on} \ \Hat{\Gamma}^\mathrm{base}\cross(0,T],
    \label{eq:bc_base_damage}
    \\
    & \bm K\nabla{\alpha}\cdot \bm N =  0, && \mbox{on} \ \Hat{\Gamma}^\mathrm{endo}\cup\Hat{\Gamma}^\mathrm{epi}\cross(0,T],
    \label{eq:bc_neu_damage}
    \end{empheq}
\end{subequations}
where Eq.~\eqref{eq:bc_endo} is the endocardial pressure boundary condition, \eqref{eq:bc_epi} models the action of the pericardial sac~\cite{pfaller_importance_2019}, and \eqref{eq:bc_base} is the energy-consistent boundary condition~\cite{regazzoni_machine_2020}. On the basal boundary $\Hat{\Gamma}^\mathrm{base}$, far away from the infarcted area, we impose a homogeneous Dirichlet boundary condition for the damage~\eqref{eq:bc_base_damage}, while on the endocardial and epicardial surfaces we impose a homogeneous Neumann boundary condition~\eqref{eq:bc_neu_damage}. We model the cardiac infarction as a local loss of active contractility $T_\mathrm{a}$ in a portion of the myocardium, according to Eq.~\eqref{eq:ta_ischemia}. The infarcted portion of myocardium is depicted in Figure~\ref{fig:infarcted_portion}. The fibre directions are determined using the Bayer et al. rule-based algorithm\footnote{The Bayer et al. rule-based algorithm \cite{bayer_novel_2012} is performed in a mesh which does not correspond to the mesh of the reference configuration $\Hat{\Omega}$, but rather an intermediate loaded imaging configuration $\Tilde{\Omega}$. Let $\Tilde{\bm F}$ be the deformation gradient of the deformation mapping $\hat{\Omega}$ to $\Tilde{\Omega}$. If the orthonormal triplet $\left\{\Tilde{\bm f}, \Tilde{\bm s}, \Tilde{\bm n}\right\}$ is obtained using the Bayer et al. algorithm on $\Tilde{\Omega}$, then the fibre directions in the reference configuration:
\begin{equation*}
    \hat{\bm f} = \frac{\Tilde{\bm F}^{-1}\Tilde{\bm f}}{\|\Tilde{\bm F}^{-1}\Tilde{\bm f}\|}, \quad \hat{\bm s} = \frac{\Tilde{\bm F}^{-1}\Tilde{\bm s}}{\|\Tilde{\bm F}^{-1}\Tilde{\bm s}\|}, \quad \hat{\bm n} = \frac{\Tilde{\bm F}^{-1}\Tilde{\bm n}}{\|\Tilde{\bm F}^{-1}\Tilde{\bm n}\|}
\end{equation*}
are no longer necessarily orthogonal, making the material no longer orthotropic, but rather generally anisotropic. In our opinion, however, this does not alter significantly the results of the simulations.} \cite{bayer_novel_2012}, and reported in Figure~\ref{fig:fibres_directions}.

\begin{figure}[t]
     \centering
     \subfloat[][]{\includegraphics[width = 0.21\textwidth]{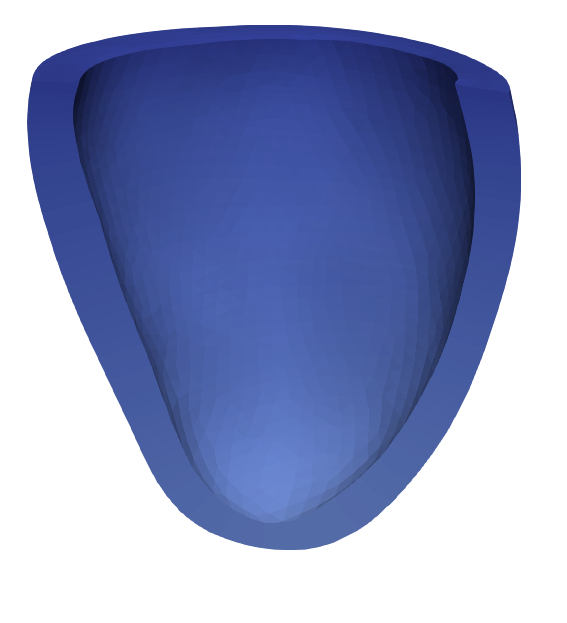}\label{fig:t_080}}\hfill
     \subfloat[][]{\includegraphics[width = 0.21\textwidth]{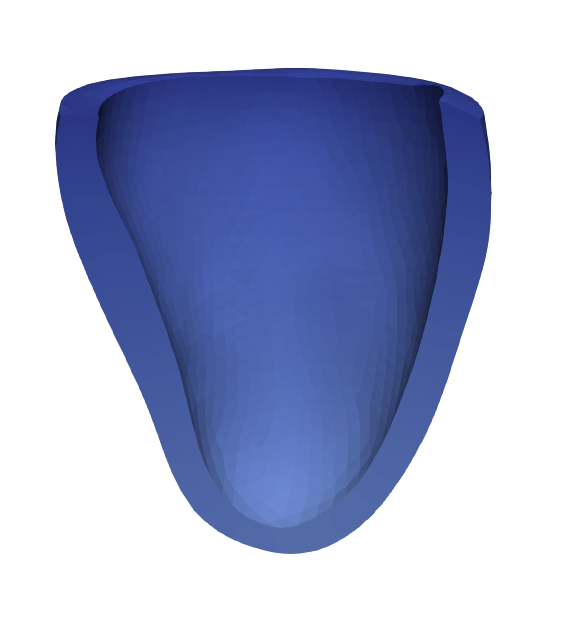}\label{fig:t_100}}\hfill
     \subfloat[][]{\includegraphics[width = 0.21\textwidth]{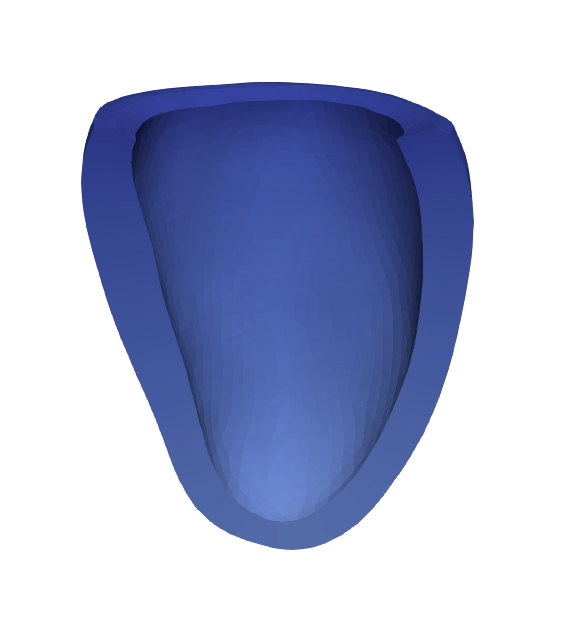}\label{fig:t_120}}\hfill
     \subfloat[][]{\includegraphics[width = 0.21\textwidth]{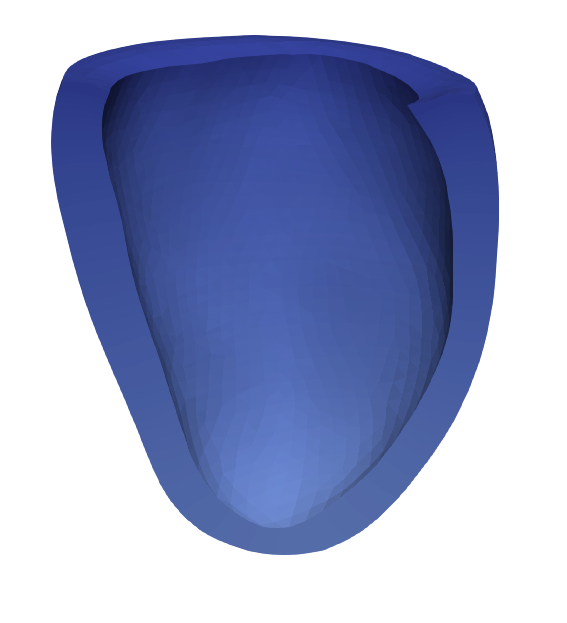}\label{fig:t_140}}\\
     \subfloat[][]{\includegraphics[width = 0.21\textwidth]{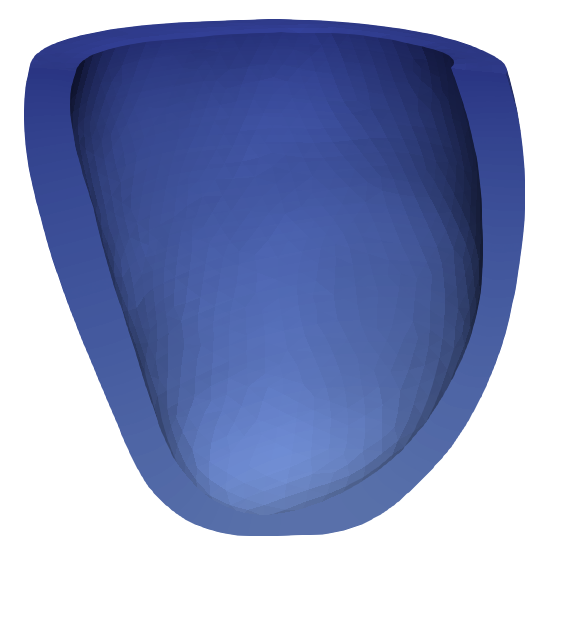}\label{fig:t_160}}\hfill
     \subfloat[][]{\includegraphics[width = 0.21\textwidth]{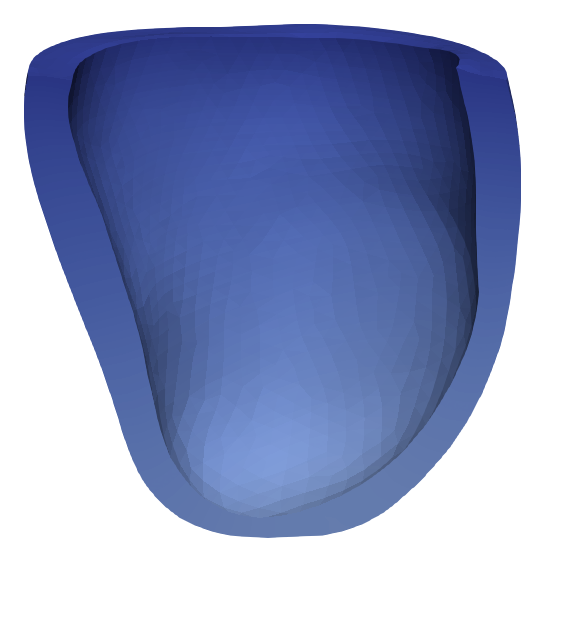}\label{fig:t_170}}\hfill
     \subfloat[][]{\includegraphics[width = 0.21\textwidth]{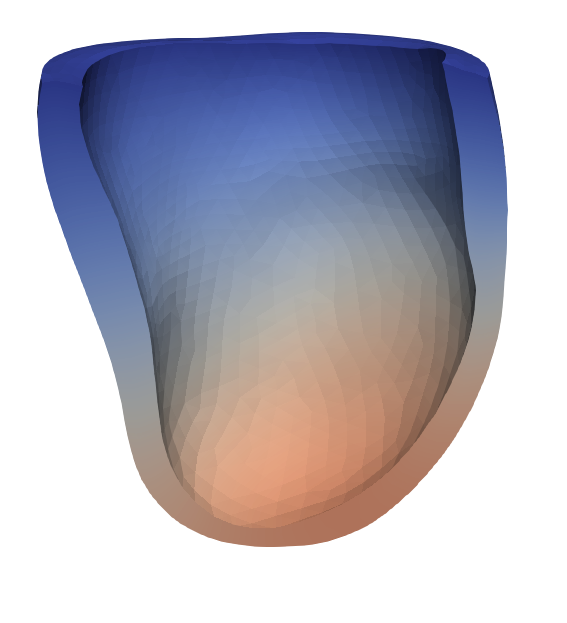}\label{fig:t_175}}\hfill
     \subfloat[][]{\includegraphics[width = 0.21\textwidth]{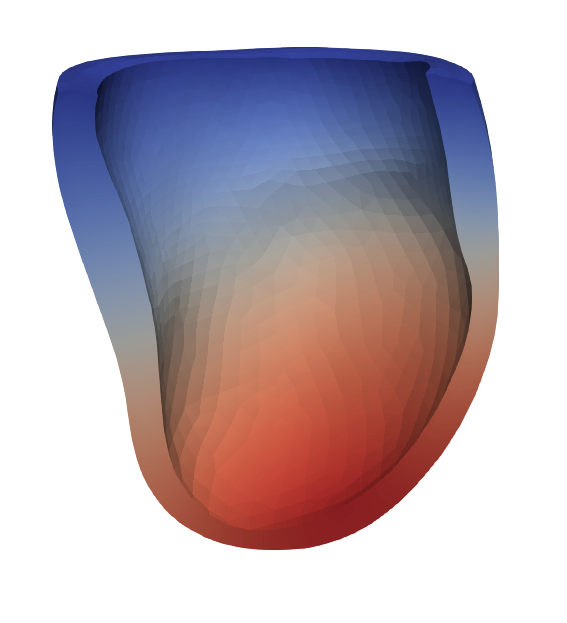}\label{fig:t_180}}\\
     \subfloat[][]{\includegraphics[width = 0.21\textwidth]{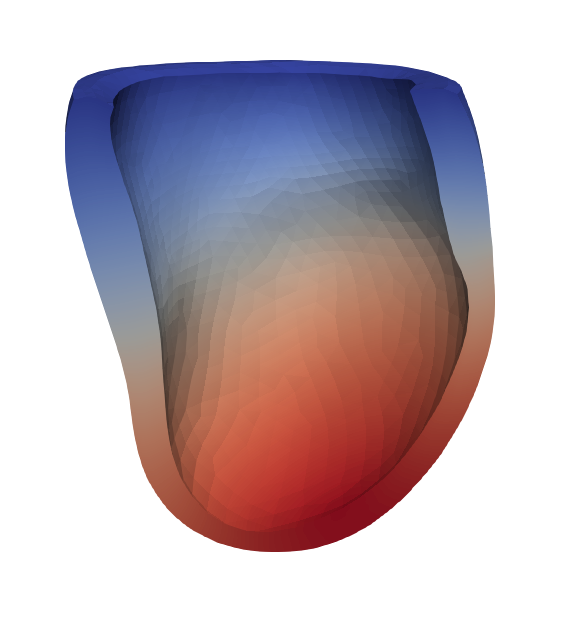}\label{fig:t_185}}\hfill
     \subfloat[][]{\includegraphics[width = 0.21\textwidth]{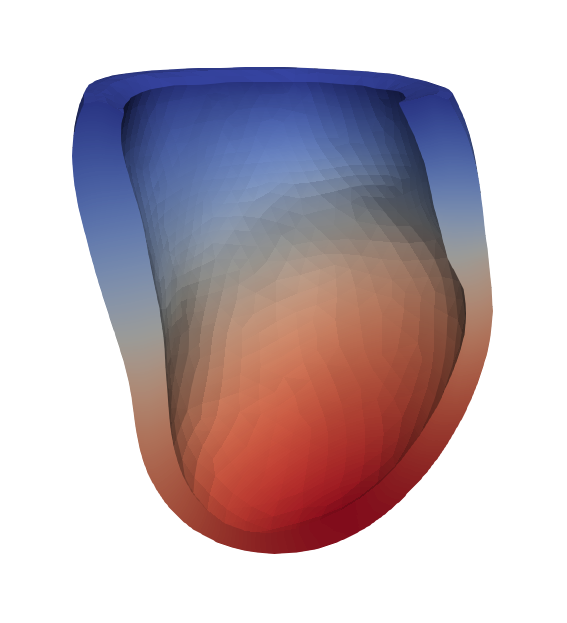}\label{fig:t_190}}\hfill
     \subfloat[][]{\includegraphics[width = 0.21\textwidth]{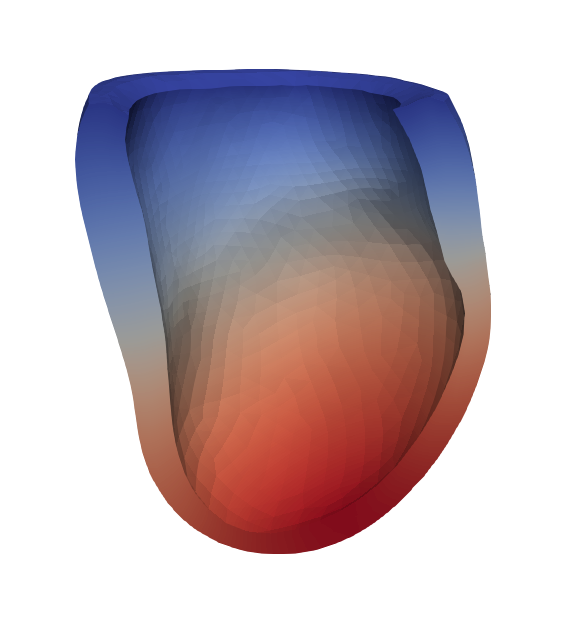}\label{fig:t_195}}\hfill
     \subfloat[][]{\includegraphics[width = 0.21\textwidth]{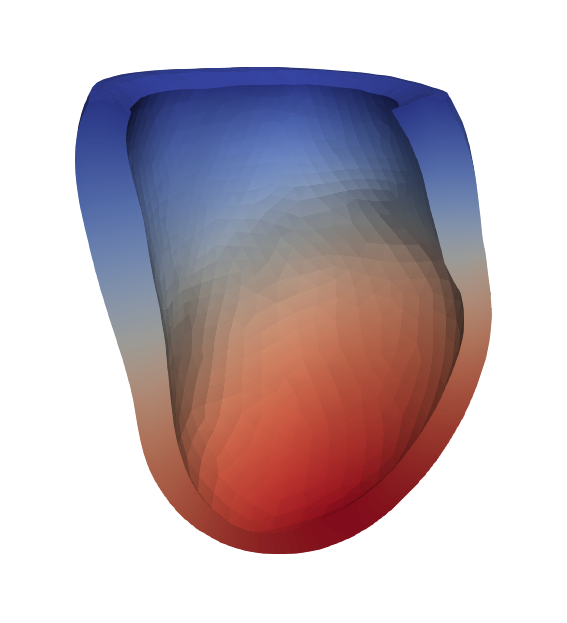}\label{fig:t_200}}\\
     \subfloat[][]{\includegraphics[width = 0.21\textwidth]{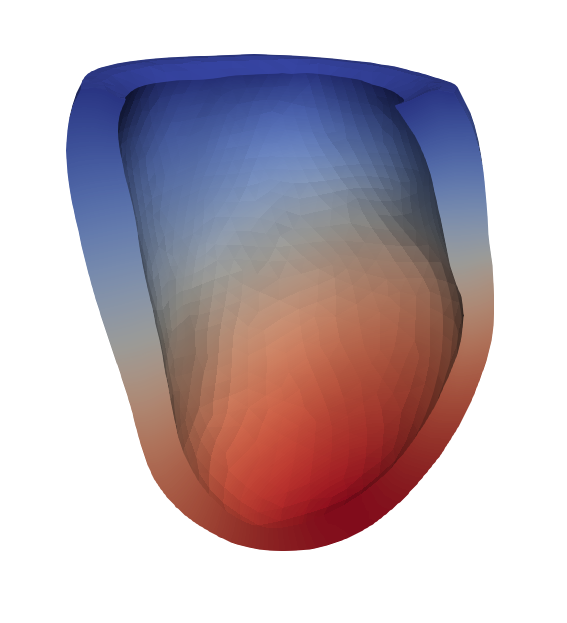}\label{fig:t_210}}\hfill
     \subfloat[][]{\includegraphics[width = 0.21\textwidth]{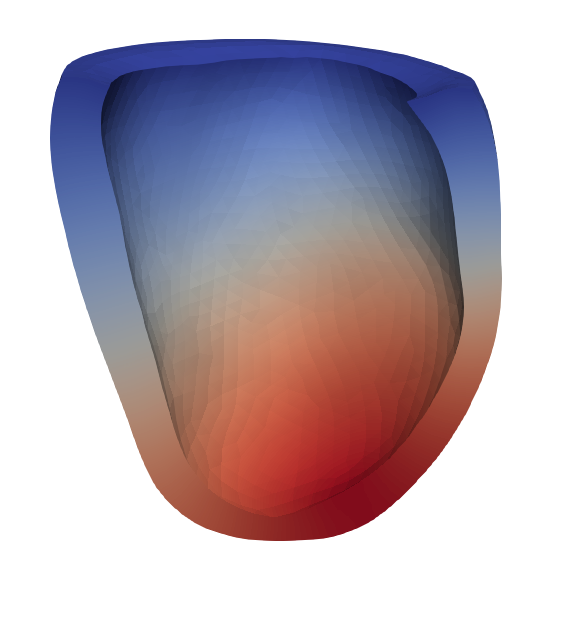}\label{fig:t_220}}\hfill
     \subfloat[][]{\includegraphics[width = 0.21\textwidth]{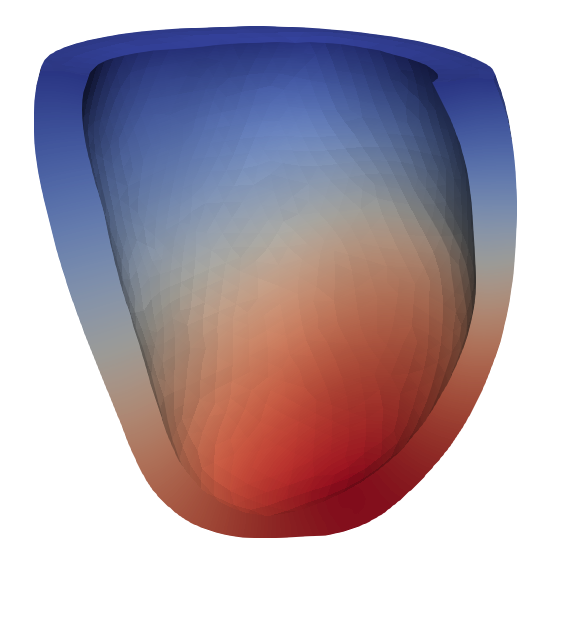}\label{fig:t_230}}\hfill
     \subfloat[][]{\includegraphics[width = 0.21\textwidth]{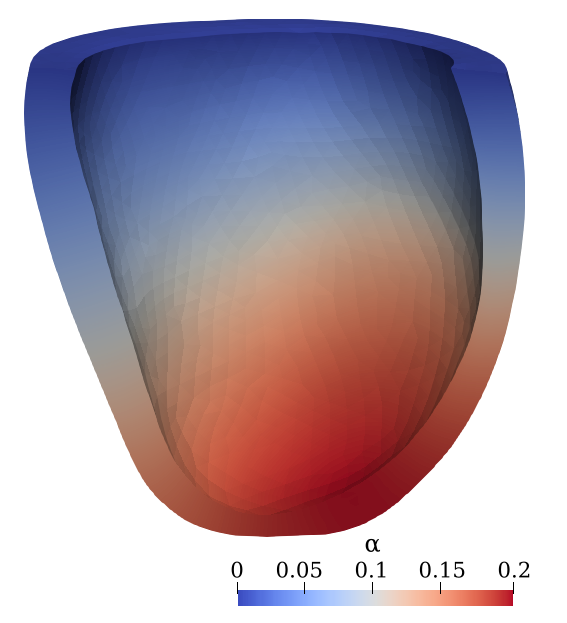}\label{fig:t_240}}\\
     \caption{Test 2: Damage field $\alpha$ evolution. (a) $t=0.80~\si{\second}$. (b) $t=1.00~\si{\second}$. (c) $t=1.20~\si{\second}$. (d) $t=1.40~\si{\second}$. (e) $t=1.60~\si{\second}$. (f) $t=1.70~\si{\second}$. (g) $t=1.75~\si{\second}$. (h) $t=1.80~\si{\second}$. (i) $t=1.85~\si{\second}$. (j) $t=1.90~\si{\second}$. (k) $t=1.95~\si{\second}$. (l) $t=2.00~\si{\second}$. (m) $t=2.10~\si{\second}$. (n) $t=2.20~\si{\second}$. (o) $t=2.30~\si{\second}$. (p) $t=2.40~\si{\second}$.}
     \label{fig:damage_evolution_lv}
\end{figure}

Following the remarks of Section~\ref{sec:residual_stresses_strains}, we model the specific fracture energy $w_1$ as space-dependent using Eq.~\eqref{eq:specific_fracture_energy_space_dep}. Details of the full left ventricle damage-electromechanics model are given in~\ref{appendix:left_ventricle_electromechanics}. Model and simulation parameters are reported in Tables~\ref{table:strain_energy_densities_params} and~\ref{table:damage_params}. We report selected timesteps for the damage evolution in the left ventricle following a myocardial infarction in Figure~\ref{fig:damage_evolution_lv}.

From Figure~\ref{fig:t_080}, we see that while the damage field is nonzero for the healthy left ventricle, it is significantly lower than following a myocardial infarction. Furthermore, Figure~\ref{fig:t_100}-\ref{fig:t_160} show that, while the there is some active function loss, since $t>t_\mathrm{init}$, the changes in internal load distribution are not sufficient to significantly affect the damage field. Figures~\ref{fig:t_170}-\ref{fig:t_240} depict the onset and evolution of the damage field in the left ventricle. Indeed, by $t=1.6~\si{\second}$ all active function has been lost in the affected region, and the affected region bulges and exhibits a systolic stretch. This stretch is highest at $t=1.9~\si{\second}$, depicted in Figure~\ref{fig:t_190}, during which the highest damage field values are attained. After the systolic stretch, in the diastole, the whole left ventricle relaxes, and damage does not appear to evolve, as evidenced by Figures~\ref{fig:t_210}-\ref{fig:t_240}. Moreover, consistently with what was experimentally observed~\cite{whittaker_role_1991}, the collagen damage is concentrated in the infarcted area.

\subsubsection{Impact on post myocardial-infarction left ventricular function}
\label{sec:results_impact}
\begin{figure}[t]
     \centering
     \subfloat[][]{\includegraphics[width = 0.48\textwidth]{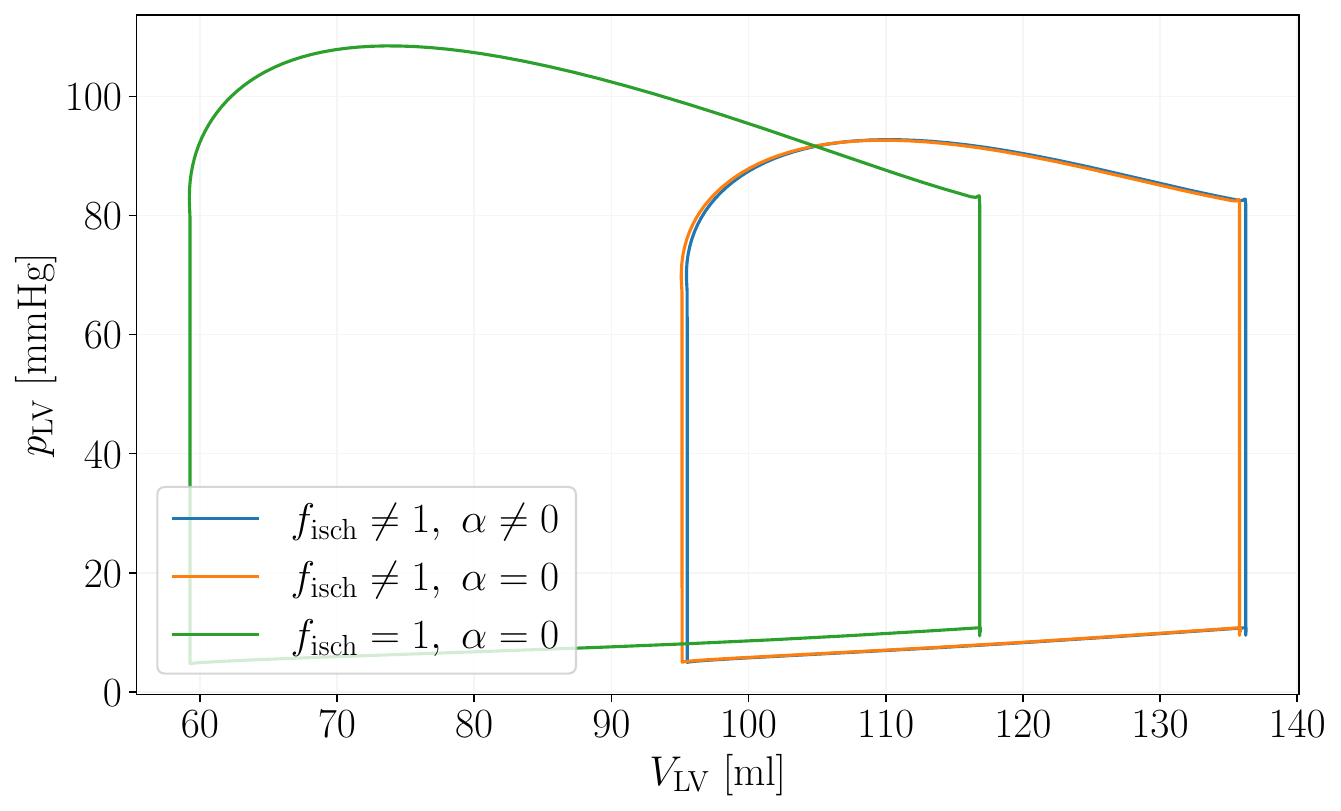}\label{fig:damage_vs_healthy}}
     \hfill
     \subfloat[][]{\includegraphics[width = 0.48\textwidth]{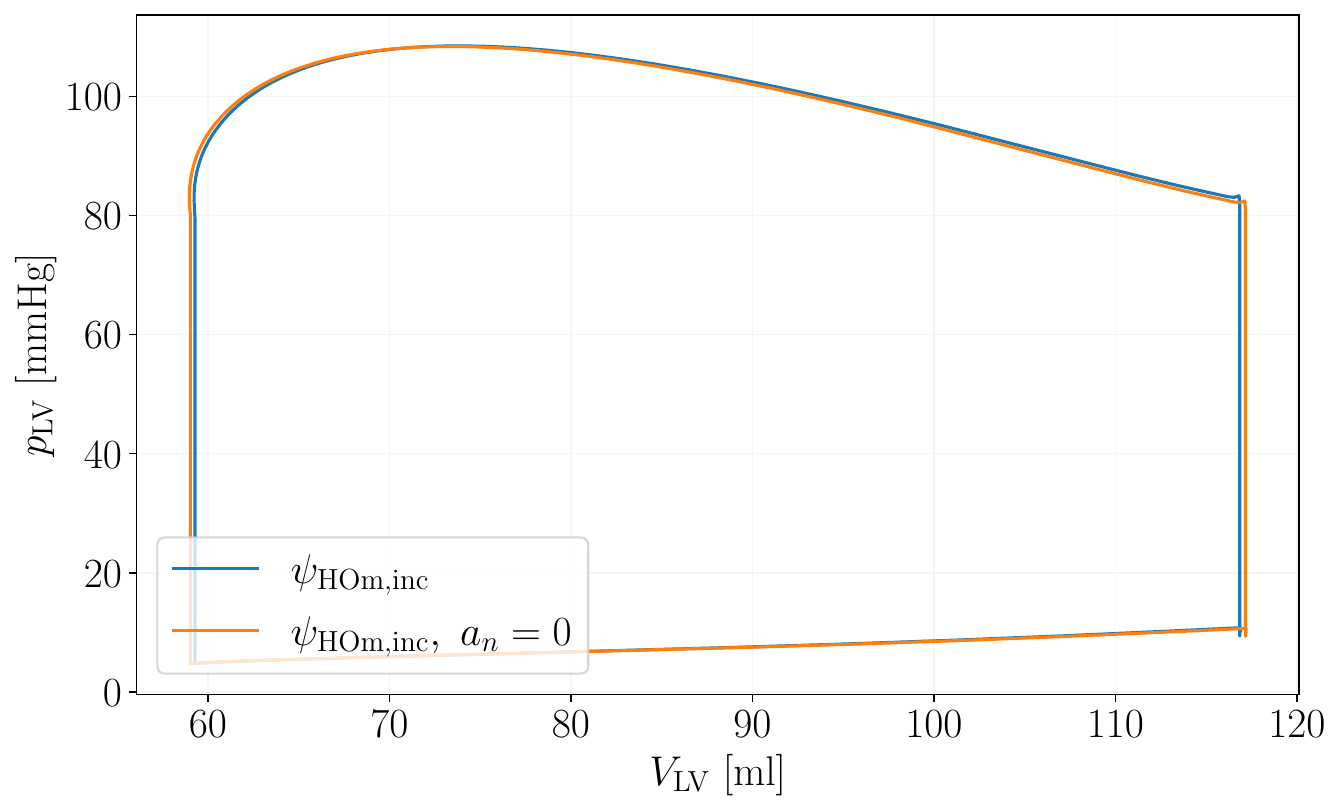}\label{fig:healthy_vs_healthy}}\\
     \caption{Test 2: Pressure-volume loops for the healthy, infarcted, and both infarcted and damaged left ventricle. (a) Comparison of healthy left ventricle, corresponding to $\alpha= 0$ and $f_\mathrm{isch}=1$, with the left ventricle where only loss of contractile function occurs, corresponding to $\alpha+ 0$ and $f_\mathrm{isch}\ne 1$, and with the left ventricle where loss of contractile function along with damage occur, corresponding to $\alpha\ne 0$ and $f_\mathrm{isch}\ne 1$. (b) Comparison of left ventricular pressure-volume loops between healthy left ventricle with $\psi_\mathrm{HOm,inc}$ strain energy density with parameters in Table~\ref{table:strain_energy_densities_params} and with the parameter $a_n=0$.}
     \label{fig:pv_loops}
\end{figure}

We study the impact of the damage on the left ventricular function in terms of pressure-volume loops and ejection fraction (EF). In Figure~\ref{fig:damage_vs_healthy}, we report the pressure-volume loops obtained for a healthy left ventricle, versus the obtained for a left ventricle which has lost local active contractility in the infarcted region, and the one in which both loss of contractility and damage occur.
From Figure~\ref{fig:damage_vs_healthy}, we see that in both cases where contractility is lost, there is 
 a rightward shift of the pressure-volume loops, and a significant worsening of the cardiac function, as evidenced by the EF reported in Table~\ref{table:pv_loops_table}. Moreover, there is an increase of the end-diastolic volume of the damaged left ventricle, owing to the softening of the tissue due to the damage. Overall, the damage appears to have little impact on the shape of the pressure-volume loop. The small effect of the damage on the pressure-volume loop is explained by Figure~\ref{fig:healthy_vs_healthy}, where we compare the pressure-volume loops of two healthy left ventricles, one of which has the material parameter in the strain-energy density $a_n$ set to zero. The differences appear to be minimal. However, Table~\ref{table:pv_loops_table} shows that in the case of the healthy left ventricle, loss of stiffness in the $\hat{\bm n}$ direction is associated with slightly improved cardiac function. This is in contrast with the case of the infarcted left ventricle, where although the end diastolic volume is decreased, there is no improvement of the EF.
\begin{table}[t]
    \centering 
    \begin{tabular}{ p{11em} p{5em} p{5em} p{4em}}
     & $V_\mathrm{ED}$ & $V_\mathrm{ES}$ & $\mathrm{EF}$  \\ 
    \hline
    $\alpha= 0,\ f_\mathrm{isch}=1$ & $116.854~\si{\milli\liter}$ & $59.287~\si{\milli\liter}$ & $49.26~\%$ \\
    $\alpha= 0,\ f_\mathrm{isch}\ne 1$ & $135.797~\si{\milli\liter}$ & $95.153~\si{\milli\liter}$ & $29.93~\%$ \\
    $\alpha\ne 0,\ f_\mathrm{isch}\ne 1$ & $136.250~\si{\milli\liter}$ & $95.529~\si{\milli\liter}$ & $29.89~\%$ \\
    $\alpha= 0,\ f_\mathrm{isch}=1,\ a_n=0$ & $117.169~\si{\milli\liter}$ & $59.033~\si{\milli\liter}$ & $49.62~\%$  \\
    \end{tabular}
    \\[10pt]
    \caption{End-diastolic volume ($V_\mathrm{ED}$), end-systolic volume ($V_\mathrm{Es}$) and ejection fraction ($\mathrm{EF}$) for the pressure-volume loops reported in Figure~\ref{fig:pv_loops}, computed as $EF = (V_\mathrm{ED}-V_\mathrm{ES})/V_\mathrm{ED}$.}
    \label{table:pv_loops_table}
\end{table}
\afterpage{\clearpage}
\section{Conclusions and discussion}
\label{sec:discussion}
In this work, we derived a model for the irreversible evolution of damage within a compressible hyperelastic medium and applied it to the damage of collagen fibres in myocardial tissue. By simulating it in a slab of tissue which was subjected to an indenter-like loading, we found that damage magnitude and evolution depend on the fibre orientation in the sample. Then, we simulated the damage of the collagen fibres within cleavage planes following an acute myocardial infarction. We have found that damage could evolve in the infarcted zone and that its evolution occurred principally during the systole, due to the increased load in the systolic stretch and the bulging out of the infarcted area.
\par
Regarding the effect of the damage, while in the healthy myocardium loss of stiffness in the cleavage planes is associated with a slightly improved cardiac function and output, the loss of stiffness in the cleavage planes which occurs in the infarcted myocardium leads to a slight worsening of the cardiac function.
\subsection{Limitations}
\label{sec:limitations}
Due to the anisotropic nature of the myocardium, we expect that the damage is also anisotropic. While this study incorporates the anisotropic nature of the damage through the dissipation energy function and the anisotropic length scale diffusion tensor $\bm K$ in the total mechanical energy, true damage anisotropy can be described only using a tensorial damage variable.
\par
Regarding the evolution of the damage field in both the case of the loaded slab and the infarcted left ventricle, the damage field remains overall low. The damage field evolution is sensitive to the choice of dissipation energy. Indeed, the assumption that only collagen fibres in the cleavage planes are affected by damage is restrictive. Higher loadings would be needed to further increase the damage field in the left ventricle. This could be achieved by increasing the contractility of the healthy tissue following the infarction by incorporating baroreflex regulation as in~\cite{sharifi_multiscale_2024}. Other than determining damage magnitude, a different choice of dissipation energy could alter the results of Section~\ref{sec:results_impact} and have a detrimental effect on cardiac function.
\par
As discussed in Section~\ref{sec:residual_stresses_strains}, the myocardial tissue is residually stressed, affecting the distribution and magnitude of the strain energy density during the cardiac cycle. The inclusion of residual stresses would yield a more accurate energy distribution and possibly allow for a single homogenous specific fracture energy to be used for the whole myocardium.
\subsection{Further developments}
\label{sec:further_developments}
As our results suggest, material degradation may not be sufficient to imply deterioration of cardiac function after a myocardial infarction other than what is due to the loss of local contractility. On the other hand, the inclusion of an inelastic deformation as described in Section~\ref{sec:introduction} may have a larger effect. Future studies will focus on the modelling of the inelastic slippage of sheetlets and whether or not it may be included in the further deterioration of cardiac function in the early stages following a myocardial infarction.
\par
As mentioned in Section~\ref{sec:limitations}, the damage field evolution depends on the choice of dissipation energy, and the assumption that only collagen fibres in the cleavage planes can be damaged is restrictive. Further studies will explore the effect of different dissipation functions, possibly containing isotropic parts or invariants of other directions, capturing both the natural dispersion of the collagen fibres' directions but also damage to different collagen fibre families contributing to the overall stiffness of the myocardium.

\section*{Acknowledgements}
IR, FR and LD acknowledge their membership to INdAM GNCS -- Gruppo Nazionale per il Calcolo Scientifico (National Group for Scientific Computing), Italy. The research of IR, FR and LD is part of the activities of ``Dipartimento di Eccellenza 2023–2027", MUR, Italy, Dipartimento di Matematica, Politecnico di Milano. 

\appendix 
\section{Boundedness of the damage variable $\alpha$}
\label{appendix:boundedness}
We prove that by replacing $g:[0,1]\to[0,1]$ with $\Tilde{g}:\mathbb R\to[0,1]$ in Eq.~\eqref{eq:energy_gradient_damage}, which is defined as:
\begin{equation}
\label{eq:truncation}
\begin{aligned}
    \Tilde{g}(\alpha) = \begin{dcases} g(1),& \mbox{if }\alpha\ge 1,\\
    g(\alpha),& \mbox{if }1>\alpha\ge 0,\\ 
    g(0), & \mbox {if } \alpha<0.\end{dcases}
\end{aligned}
\end{equation}
then,
\begin{equation*}
    \forall(\bm d,\alpha)\in V\times W\ \exists \bar{\alpha}\in W,\ \Bar{\alpha}\in[0,1]\ \mbox{a.e. in }\Hat{\Omega}\ \mbox{s.t. }\mathcal{E}({\bm d},\bar{\alpha})\leq \mathcal{E}({\bm d},{\alpha}),
\end{equation*}
Let us define:
\begin{equation*}
    \bar\alpha=\begin{dcases}1, & \mbox {in } \left\{\alpha\geq 1\right\}, \\ \alpha, & \mbox {in } \left\{1>\alpha\geq 0\right\}, \\ 0, & \mbox {in } \left\{\alpha<0\right\}.\end{dcases}
\end{equation*}
Then,
\begin{equation*}
\begin{aligned}
    \mathcal{E}({\bm d},{\alpha}) - \mathcal{E}({\bm d},\bar{\alpha}) &= \int_{\left\{\alpha>1\right\}} \psi_\mathrm{diss}(\bm F)(g(1)-g(\alpha))+ \int_{\left\{\alpha>1\right\}} w_1\left((\alpha^2-1)+\ell^2|\sqrt{\bm K}\nabla\alpha|^2\right) \\&+ \int_{\left\{\alpha<0\right\}} \psi_\mathrm{diss}(\bm F)(g(0)-g(\alpha))+ \int_{\left\{\alpha<0\right\}} w_1\left(\alpha^2+\ell^2|\sqrt{\bm K}\nabla\alpha|^2\right) \\&=\int_{\left\{\alpha>1\right\}} w_1\left((\alpha^2-1)+\ell^2|\sqrt{\bm K}\nabla\alpha|^2\right) + \int_{\left\{\alpha<0\right\}} w_1\left(\alpha^2+\ell^2|\sqrt{\bm K}\nabla\alpha|^2\right),
\end{aligned}
\end{equation*}
which is strictly positive if the Lebesgue measure of $\left\{\alpha>1\right\}$ or of $\left\{\alpha<0\right\}$ is strictly greater than zero. This means that if we use the degradation function defined as in Eq.~\eqref{eq:truncation}, the minimizers automatically satisfy the bound $\alpha\in[0,1]$.

\section{On the choice of strain energy density}
\label{appendix:monotonicity_of_the_stress}
\begin{figure}[t]
     \centering
     \subfloat[][]{\includegraphics[width = 0.48\textwidth]{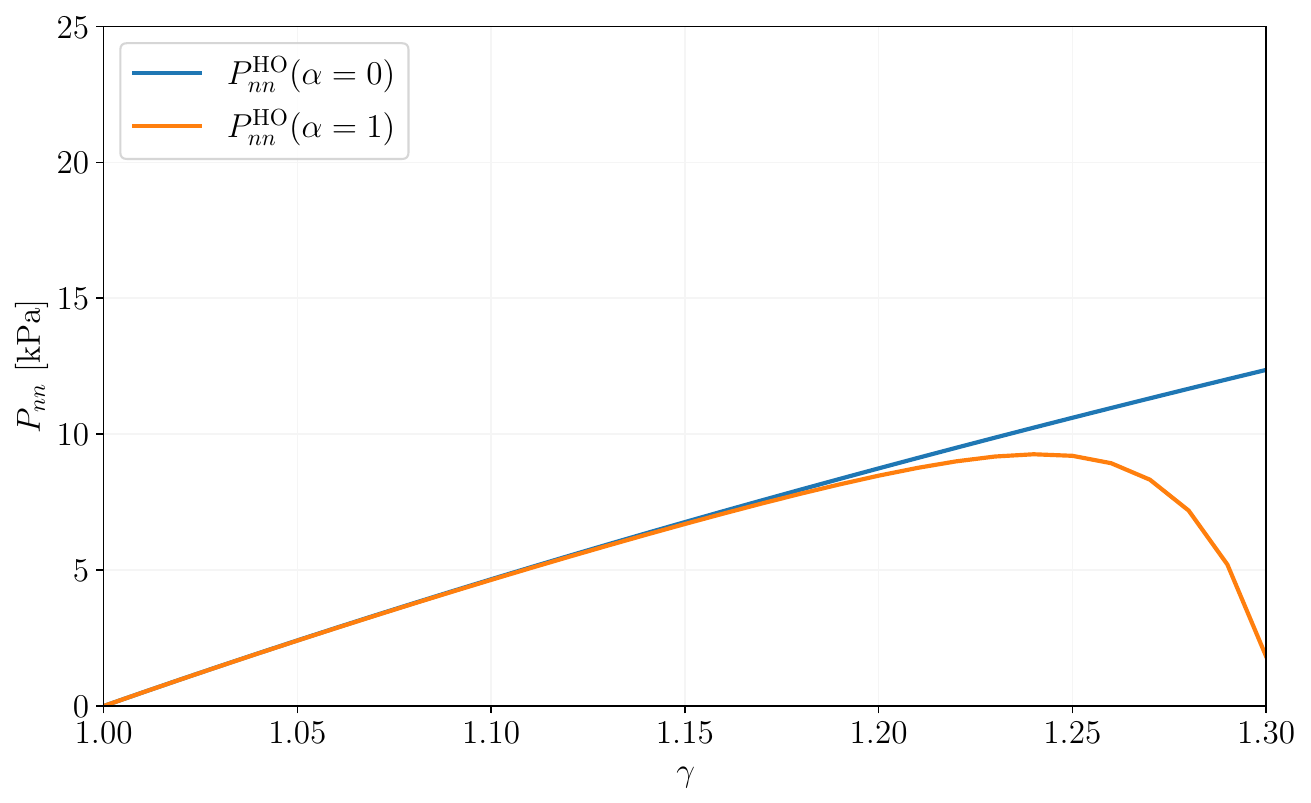}\label{fig:nonmonotonicity_nmod}}\hfill
     \subfloat[][]{\includegraphics[width = 0.48\textwidth]{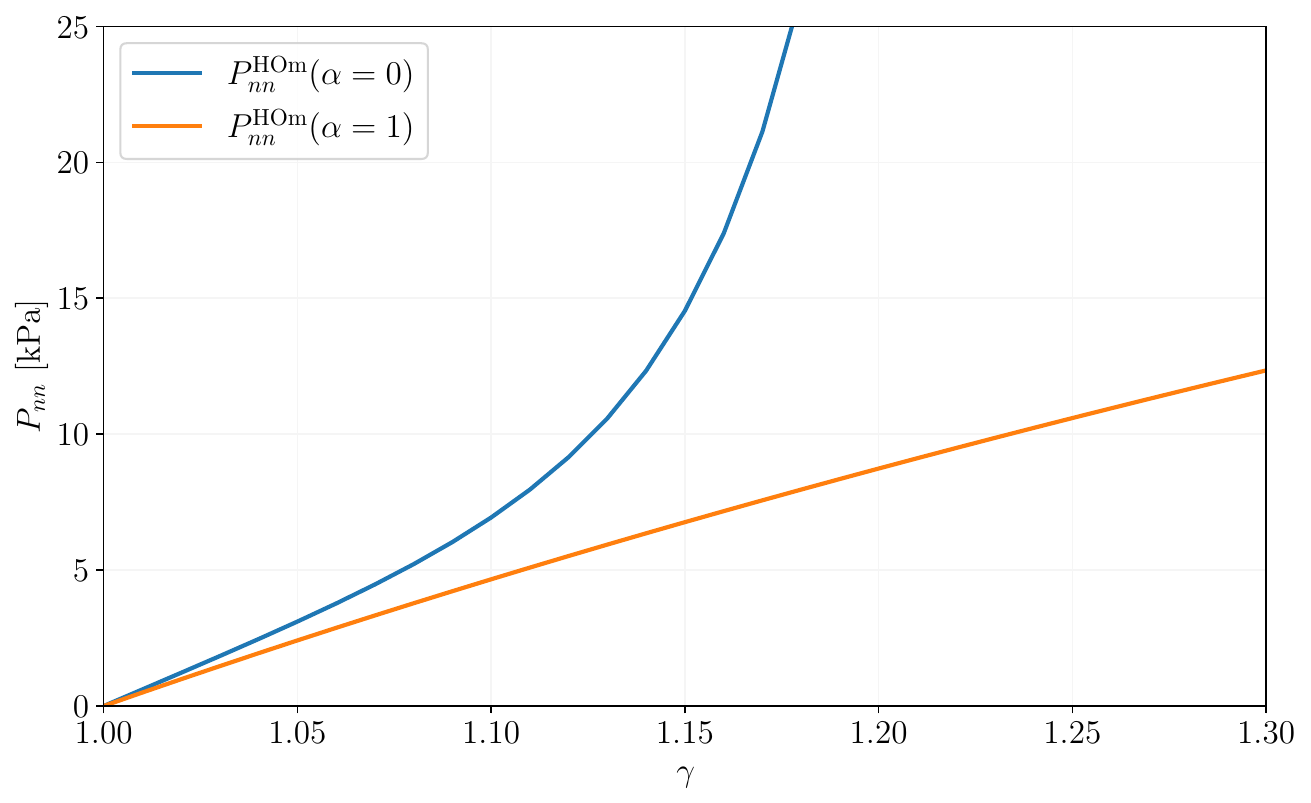}\label{fig:nonmonotonicity_mod}}\\
     \caption{$P_{nn}$ depending on the strain $\gamma$ and the damage parameter $\alpha$. (a) $\psi_\mathrm{HO}$ with isochoric-volumetric split. (b) $\psi_\mathrm{HO,inc}$. $P_{nn}^\mathrm{HOm}(\alpha=0)$ was cut for visualization purposes.}
     \label{fig:stress_monotonicity}
\end{figure}

We explain the motivation behind the modification of the existing Holzapfel-Ogden model as in~\eqref{eq:ho_mod_qinc}. For a nonlinear elastic material, the Strong Ellipticity Condition~\cite{antman_nonlinear_2005}:
\begin{equation}
    \label{eq:strong_ellipticity_condition}
    [\bm P(\bm G + \alpha\bm H) - \bm P (\bm G)] :\bm H > 0, \quad \forall \bm  G\in \mathrm{Lin}^+,\quad \forall \bm H \mbox{ of rank 1},\quad \forall\alpha\in(0, 1] \mbox{ s.t. } \mathrm{det}(\bm G + \alpha\bm H) > 0,
\end{equation}
is a constitutive assumption guaranteeing the existence of real solutions to the incremental problem~\cite{antman_nonlinear_2005}. We recall that a corollary of the Strong Ellipticity Condition for a hyperelastic material is that, under suitable differentiability assumptions for the strain energy density, for the first Piola-Kirchhoff stress tensor $\bm P$ it holds that~\cite{antman_nonlinear_2005}:
\begin{equation}
\label{eq:strong_ellipticity_corollary}
    \pdv{(\bm a\cdot \bm P(\bm F) \bm b)}{(\bm a \cdot \bm F \bm b)}>0,
\end{equation}
for any $\bm a $ and $\bm b$ that are independent of $\bm a\cdot\bm F\bm b$.
For a material with purely elastic part given by~\eqref{eq:ho_mod_qinc}, the Piola-Kirchhoff stress tensor is given by:
\begin{equation*}
    \bm P = \bm P_1 + \bm P_f + \bm P_s + \bm P_n + \bm P_{fs} +  \bm P_J,
\end{equation*}
where
\begin{equation*}
    \begin{aligned}
        &\bm P_1 = a\exp\left(b(\Bar{I}_1-3)\right)\left(\frac{1}{J^{2/3}}\bm F-\frac{1}{3}\Bar{I}_1\bm F^{-T}\right), \\
        &\bm P_f = {2}a_f\langle I_{4f}-1\rangle_+\exp\left(b_f\langle I_{4f}-1\rangle_+^2\right)\bm F\Hat{\bm f}\otimes\Hat{\bm f},\\
        &\bm P_s = {2}a_s\langle I_{4s}-1\rangle_+\exp\left(b_s\langle I_{4s}-1\rangle_+^2\right)\bm F\Hat{\bm s}\otimes \Hat{\bm s},\\
        &\bm P_n = 2(1-\alpha)^2a_n\langle I_{4n}-1\rangle_+\exp\left(b_n\langle I_{4n}-1\rangle_+^2\right)\bm F\Hat{\bm n}\otimes \Hat{\bm n},\\
        &\bm P_{fs} = a_{fs}I_{8fs}\exp\left(b_{fs}I_{8fs}^2\right)\left(\bm F\Hat{\bm f}\otimes \Hat{\bm s} + \bm F\Hat{\bm s}\otimes \Hat{\bm f}\right),\\
        &\bm P_{J} = \frac{c_\mathrm{bulk}}{2}\left(J\log(J)+J-1\right)\bm F^{-T}.
    \end{aligned}
\end{equation*}
On the other hand, the original Holzapfel-Ogden model~\cite{holzapfel_constitutive_2009} without explicit stiffening in the $\hat{\bm n}$ direction with the dissipation in the cleavage planes given by Eq.~\eqref{eq:psi_diss_def} yields the stress $\bm P^\mathrm{HO}_n$, that takes the form of:
\begin{equation*}
    \bm P^\mathrm{HO}_n = 2((1-\alpha)^2-1)a_n\langle I_{4n}-1\rangle_+\exp\left(b_f\langle I_{4n}-1\rangle_+^2\right)\bm F\Hat{\bm n}\otimes \Hat{\bm n}.
\end{equation*}
We compare how these two material models satisfy condition~\eqref{eq:strong_ellipticity_corollary}. In particular, let us choose $\bm a=\bm b = \Hat{\bm n}$ and $\bm F = \Hat{\bm f}\otimes\Hat{\bm f} + \Hat{\bm s}\otimes\Hat{\bm s} + \gamma \Hat{\bm n}\otimes\Hat{\bm n}$ with $\gamma>1$, so that the derivative appearing in~\eqref{eq:strong_ellipticity_corollary} is well-defined. Let us also assume that locally the material is fully damaged, $\alpha(\bm X) = 1$. Then, for the model without explicit stiffening in the $\hat{\bm n}$ direction, it holds that:
\begin{equation*}
\begin{aligned}
    &\pdv{(\Hat{\bm n}\cdot \bm P(\bm F) \Hat{\bm n})}{(\Hat{\bm n} \cdot \bm F \Hat{\bm n})} = \pdv{(\Hat{\bm n}\cdot \bm P(\bm F) \Hat{\bm n})}{\gamma} = \pdv{(\Hat{\bm n}\cdot \bm P_1(\bm F) \Hat{\bm n})}{\gamma}+\pdv{\Hat{\bm n}\cdot \bm P^\mathrm{HO}_n(\bm F) \Hat{\bm n}}{\gamma}+\pdv{\Hat{\bm n}\cdot \bm P_J(\bm F) \Hat{\bm n}}{\gamma}
    \\ &= \frac{a}{3}\exp\left[b(2\gamma^{-\frac{2}{3}}+\gamma^{\frac{4}{3}}-3)\right]\left(-4b\gamma^{-\frac{4}{3}}+\frac{8b}{3}\gamma^{-\frac{7}{3}}-\frac{4b}{3}\gamma^{-\frac{1}{3}}+4b\gamma^{\frac{2}{3}}-\frac{4b}{3}\gamma^{\frac{5}{3}}+\gamma^{-\frac{2}{3}}+\frac{4}{3}\gamma^{-\frac{5}{3}}-\frac{4}{3}\gamma^{\frac{1}{3}}\right)  \\ &-2a_n\exp\left[b_n(\gamma^2-1)^2\right] \left(4b_n\gamma^6-8b_n\gamma^4+(4b_n+3)\gamma^2-1\right) \\&+\frac{c_\mathrm{bulk}}{2}\left(\gamma^{-1}+\gamma^{-2}\right).
\end{aligned}
\end{equation*}
In particular, the last term, coming from $\bm P_J$ is always positive and bounded, while the second term, coming from $\bm P_n$, is negative for $\gamma$ large enough for any $a_n>0$ and $b_n>1$. On the other hand, for material parameters such that $b\centernot\gg b_n$ and $a\centernot\gg a_n$, that is, where the matrix is not significantly exponentially stiffer than the fibres, the first term, coming from $\bm P_1$, becomes negative for $\gamma$ large enough. Moreover, the contribution of the second term increases significantly faster in $\gamma$, making so that it has adverse effects on the range of strains $\gamma$ over which the Strong Ellipticity Condition has a chance to hold. On the other hand, for the constitutive model incorporating the explicit stiffening in the direction $\hat{\bm n}$ and full damage, condition~\eqref{eq:strong_ellipticity_corollary} becomes: 
\begin{equation*}
\begin{aligned}
    &\pdv{(\Hat{\bm n}\cdot \bm P(\bm F) \Hat{\bm n})}{(\Hat{\bm n} \cdot \bm F \Hat{\bm n})} = \pdv{(\Hat{\bm n}\cdot \bm P(\bm F) \Hat{\bm n})}{\gamma} = \pdv{(\Hat{\bm n}\cdot \bm P_1(\bm F) \Hat{\bm n})}{\gamma}+\pdv{(\Hat{\bm n}\cdot \bm P^\mathrm{HO}_n(\bm F) \Hat{\bm n})}{\gamma}+\pdv{(\Hat{\bm n}\cdot \bm P_J(\bm F) \Hat{\bm n})}{\gamma}
    \\ &= \frac{a}{3}\exp\left[b(2\gamma^{-\frac{2}{3}}+\gamma^{\frac{4}{3}}-3)\right]\left(-4b\gamma^{-\frac{4}{3}}+\frac{8b}{3}\gamma^{-\frac{7}{3}}-\frac{4b}{3}\gamma^{-\frac{1}{3}}+4b\gamma^{\frac{2}{3}}-\frac{4b}{3}\gamma^{\frac{5}{3}}+\gamma^{-\frac{2}{3}}+\frac{4}{3}\gamma^{-\frac{5}{3}}-\frac{4}{3}\gamma^{\frac{1}{3}}\right)  \\&+\frac{c_\mathrm{bulk}}{2}\left(\gamma^{-1}+\gamma^{-2}\right).
\end{aligned}
\end{equation*}
The modified strain-energy density~\eqref{eq:ho_mod_qinc} behaves better at smaller strains when coupled with damage to the fibres the $\Hat{\bm n}$ direction, possibly guaranteeing existence pf real solutions to the incremental problem over a larger range of strains. The lack of monotonicity of the Piola-Kirchhoff stress with respect to the corresponding strain can also be shown numerically for a set of given parameters. In particular, in Figure~\ref{fig:stress_monotonicity} we report the values of the stress $P_{nn}=\hat{\bm n}\cdot\bm P\hat{\bm n}$ depending on the value of the strain $\gamma$.
In particular, we compute the values of $P_{nn}$ in two cases. First, in Figure~\ref{fig:nonmonotonicity_nmod}, we consider the original Holzapfel-Ogden strain energy $\psi_\mathrm{HO}$ with parameters given in Table~\ref{table:strain_energy_densities_params} and a dissipation energy given by Eq.~\eqref{eq:psi_diss_def} with material parameter $a_n = 29.5$ and $b_n=11.12$. As evidenced by~\ref{fig:nonmonotonicity_nmod}, the stress is not monotonically increasing with respect to the strain $\gamma$. On the other hand, in Figure~\ref{fig:nonmonotonicity_mod}, we depict the modified Holzapfel-Ogden model with dissipation energy~\eqref{eq:psi_diss_def}, with both sets of parameters given in Table~\ref{table:strain_energy_densities_params}. In particular, we see from Figure~\ref{fig:nonmonotonicity_mod} that no loss of monotonicity occurs.

\section{Left ventricle electromechanics model}
\label{appendix:left_ventricle_electromechanics}
We rely on the cardiac electromechanics model presented in~\cite{stella_fast_2022}. The full left ventricle damage-electromechanics model reads:
\begin{subequations}
    \label{eq:full_em_space_appx}
    \begin{empheq}[left=\empheqlbrace]{align}
    & \pdv{\boldsymbol{y}}{t} = \boldsymbol{h}\left(\boldsymbol{y},z_\mathrm{Ca}(t-\tau_\mathrm{a}),{\mathrm{SL}},{\pdv{{\mathrm{SL}}}{t}}\right), & \mbox{in} \ \Hat\Omega\cross(0,T], \label{eq:full_sarcomere_variables_space}
    \\
    & \rho\pdv[2]{\bm{d}}{t} - \nabla\cdot\bm{P}\left(\bm{d},T_\mathrm{a}\left(\boldsymbol{y},\mathrm{SL}\right)\right) = \bm 0, & \mbox{in} \ \Hat\Omega\cross(0,T],
    \label{eq:full_mechanics}
    \\
    & - g'(\alpha)\xi(t)+w_1{\psi_w'(\alpha)}- 2w_1\ell^2\nabla\cdot\left(\bm K\nabla\alpha\right) = 0, & \mbox{in} \ \Hat{\Omega}\cross(0,T],
    \label{eq:full_damage_mechanics}
    \\
    & \mathcal{C}\left(p_\mathrm{LV},V_\mathrm{LV}(\bm{d}),t\right) = 0, & \mbox{in} \ \Hat\Omega\cross(0,T],
    \label{eq:full_circulation}
    \\
    & \bm{P}\left(\bm{d},T_\mathrm{a}\left(\boldsymbol{y},\mathrm{SL}\right)\right)\bm{N} = -p_\mathrm{LV}(t)J\bm{F}^{-T}\bm{N}, & \mbox{on} \ \Hat\Gamma^\mathrm{endo}\cross(0,T], \label{eq:full_bc_endo}
    \\
    & \bm{P}\left(\bm{d},T_\mathrm{a}\left(\boldsymbol{y},\mathrm{SL}\right)\right)\bm{N} + \mathbf{K}^\mathrm{epi}\bm{d} + \mathbf{C}^\mathrm{epi}\pdv{\bm{d}}{t} = \mathbf{0}, & \mbox{on} \ \Hat\Gamma^\mathrm{epi}\cross(0,T],
    \label{eq:full_bc_epi}
    \\
    & \bm{P}\left(\bm{d},T_\mathrm{a}\left(\boldsymbol{y},\mathrm{SL}\right)\right)\bm{N} = p_\mathrm{LV}(t)||J\bm{F}^{-T}\bm{N}|| \frac{\int_{\Hat{\Gamma}^\mathrm{endo}}J\bm{F}^{-T}\bm{N} \mathrm{d}A }{\int_{\Hat\Gamma^\mathrm{base}}||J\bm{F}^{-T}\bm{N}||\mathrm{d}A}, & \mbox{on} \ \Hat\Gamma^\mathrm{base}\cross(0,T],
    \label{eq:full_bc_base}
    \\
    & \left(u,\boldsymbol{z},\boldsymbol{y},\mathbf{d},\dot{\mathbf{d}}\right)(0)=\left(u_0,\boldsymbol{z}_0,\boldsymbol{y}_0,\mathbf{d}_0,\mathbf{0}\right), & \mbox{in} \ \Hat\Omega,
    \label{eq:full_ic}
    \end{empheq}
\end{subequations}
where the activation time $\tau_\mathrm{a} = \tau_\mathrm{a}(\bm X)$ is obtained by solving the Eikonal-Diffusion model~\cite{stella_fast_2022},
\begin{subequations}
    \label{eq:eikonal}
    \begin{empheq}[left=\empheqlbrace]{align}
    & c_0\sqrt{\nabla\tau_\mathrm{a}\cdot\frac{1}{\chi_\mathrm{m} C_\mathrm{m}}\mathbf{D}\nabla\tau_\mathrm{a}} - \nabla\cdot\left(\frac{1}{\chi_\mathrm{m} C_\mathrm{m}}\mathbf{D}\nabla\tau_\mathrm{a}\right) = 1, \hspace{0.75em} &&\mbox{in} \ \Hat{\Omega},
    \label{eq:eik_system}
    \\
    & \left(\frac{1}{\chi_\mathrm{m} C_\mathrm{m}}\mathbf{D}\nabla\tau_\mathrm{a}\right)\cdot\mathbf{N} = 0, \hspace{0.75em} &&\mbox{on} \ \partial\Hat{\Omega}\setminus\partial\Hat{\Omega}_\mathrm{a}, \label{eq:eik_boundary}
    \\
    & \tau_\mathrm{a} = t_\mathrm{a}, \hspace{0.75em} &&\mbox{on} \ \partial\Hat{\Omega}_\mathrm{a},
    \label{eq:eik_dirichlet}
    \end{empheq}
\end{subequations}
and where Eq.~\eqref{eq:full_sarcomere_variables_space} is the system of ordinary differential equations governing the active force transient $\Hat T_\mathrm{a}(t)$,
\begin{equation*}
    \Hat T_\mathrm{a}(t) = \Hat T_\mathrm{a}(\bm y, \mathrm{SL}, \pdv{\mathrm{SL}}{t}),
\end{equation*}
and depends on the intracellular calcium concentration transient $z_\mathrm{Ca}$ and the current sarcomere length $\mathrm{SL}$,
\begin{equation*}
    \mathrm{SL} = \mathrm{SL}_0\sqrt{I_{4f}}.
\end{equation*}
In particular, for the active force generation given by Eq.~\eqref{eq:full_sarcomere_variables_space} we use the Regazzoni et al. 2020 mean-field model~\cite{regazzoni_biophysically_2020}.
The calcium concentration transient
\begin{equation*}
    z_\mathrm{Ca}:[0,T]\to\mathbb R,
\end{equation*}
is obtained as solution to the set of ordinary differential equations given by the ten Tusscher-Panfilov 2006~\cite{ten_tusscher_alternans_2006} ionic model:
\begin{subequations}
    \label{eq:ionic}
    \begin{empheq}[left=\empheqlbrace]{align}
    & \dv{\boldsymbol{z}}{t} = \boldsymbol{g}\left(u, \boldsymbol{z}\right), \hspace{0.75em} &&\mbox{in} \ (0,T], \label{eq:ionic_variables}
    \\
    & C_\mathrm{m}\dv{u}{t} + \mathcal{I}_\mathrm{ion}\left(u,\boldsymbol{z}\right) = \mathcal{I}_{\mathrm{app}}(t), \hspace{0.75em} &&\mbox{in} \ (0,T], \label{eq:ionic_cable}
    \\
    & (\bm z,u)(0) = (\bm z_0,u_0) \label{eq:ionic_init_cond}.
    \end{empheq}
\end{subequations}
The applied current $\mathcal{I}_{\mathrm{app}}(t)$ in Eq.~\eqref{eq:ionic} takes the form of:
\begin{equation}
    \label{eq:applied_current}
    \mathcal{I}_\mathrm{app}(t) = \Bar{\mathcal{I}}_\mathrm{app}\chi_{\left\{mT_\mathrm{hb}<t<mT_\mathrm{hb}+t_\mathrm{app} \right\}}(t), \quad t\in\left[0,T \right],\ m\in\mathbb Z,
\end{equation}
where $\chi_{\left\{a<t<b \right\}}(t)$ is the indicator function of the interval $[a,b]$, $T_\mathrm{hb}$ is the heartbeat period, $t_\mathrm{app}$ is the current duration, and $\Bar{\mathcal{I}}_\mathrm{app}>0$ is a current amplitude.
\par
The left ventricular pressure $p_\mathrm{LV}(t)$ is yielded by its relationship with the intraventricular contour volume $V_\mathrm{LV}(\bm d)$, given by the circulation model $\mathcal{C}(p_\mathrm{LV},V_\mathrm{LV}(\bm d),t)$ which depends on the phase of the cardiac cycle. The intraventricular volume is computed as from the displacement using the centerline method~\cite{rossi_anisotropic_2014}. The four phases of the cardiac cycle consist of the isovolumic contraction, the ejection phase, the isovolumic relaxation and the filling phase. For the two isovolumic phases, the circulation model is that $V_\mathrm{LV}(\mathbf{d})$ remain constant, for which the pressure $p_\mathrm{LV}$ acts as a Lagrange multiplier. The isovolumic contraction lasts until the aortic valve opening pressure $p_\mathrm{AVO}$ is reached. For the the ejection phase, we use the two-element windkessel model \cite{westerhof_arterial_2009} which relates the intraventricular pressure and volume as:
\begin{subequations}
    \label{eq:windkessel}
    \begin{empheq}[left=\empheqlbrace]{align}
    & C_\mathrm{WK}\dv{p_\mathrm{LV}}{t}=-\frac{p_\mathrm{LV}}{R_\mathrm{WK}}-\dv{V_\mathrm{LV}}{t},\quad t\in\left(t_\mathrm{AVO},t_\mathrm{AVC}\right),\\
    &p_\mathrm{LV}(t_\mathrm{AVO})=p_\mathrm{AVO},
    \end{empheq}
\end{subequations}
where $C_\mathrm{WK},R_\mathrm{WK},p_\mathrm{AVO}$ are model parameters.  {The ejection phase lasts until $\dv{V_\mathrm{LV}}{t}\approx 0$. Then, the isovolumic relaxation phase lasts until the pressure falls below the mitral valve opening pressure $p_\mathrm{MVO}$. Finally, the filling phase is modelled by imposing a linearly increasing pressure $p_{LV}$ until the initial end-diastolic pressure $p_\mathrm{LV}(0)$ is reached.}
\par
The function $f_\mathrm{isch}(\bm X,t)$ modelling loss of active function due to cardiac infarction takes the form of:
\begin{equation}
\label{eq:def_f_isch}
    f_\mathrm{isch}(\bm X,t) = \begin{dcases}1-f_\mathrm{pz}(\bm X), & \mbox {if } t \geq t_\mathrm{fin}, \\ 1-f_\mathrm{pz}(\bm X)\frac{t-t_\mathrm{init}}{t_\mathrm{fin}-t_\mathrm{init}}, & \mbox {if } t_\mathrm{fin} > t \geq t_\mathrm{init} , \\ 1, & \mbox {if } t < t_\mathrm{init}.\end{dcases}
\end{equation}
where $f_\mathrm{pz}(\bm X)$ is the process zone function, which takes the form of:
\begin{equation*}
    f_\mathrm{pz}(\bm X)={\chi_{\mathcal{B}_{r_\mathrm{pz}}(\bm X_\mathrm{pz})}}(\bm X) + {\chi_{\mathcal{B}_{r_\mathrm{pz}+d_\mathrm{pz}}(\bm X_\mathrm{pz})\setminus\mathcal{B}_{r_\mathrm{pz}}(\bm X_\mathrm{pz})}}(\bm X)\frac{r_\mathrm{pz}+d_\mathrm{pz}-d(\bm X,\bm X_\mathrm{pz})}{d_\mathrm{pz}},
\end{equation*}
and takes the value $1$ for all material points $\bm X$ within radius $r_\mathrm{pz}$ of the centre of the process zone $\bm X_\mathrm{pz}$, takes the value $0$ for all points with distance from the centre of the process zone $d(\bm X,\bm X_\mathrm{pz})$ larger than $d_\mathrm{pz}$, and takes intermediate values in the remaining border zone of thickness $d_\mathrm{pz}$.

\bibliographystyle{elsarticle-num} 
\bibliography{cas-refs}
\biboptions{numbers,sort&compress}




\end{document}